\documentclass[11pt,twoside]{article}  
\usepackage{hyperref}
\usepackage[ansinew]{inputenc}
\usepackage{amsfonts, color}
\usepackage{amssymb,amsmath}
\usepackage{latexsym}
\usepackage{graphicx}

\usepackage{bm}

\def\E{\mathbb{E}}
\def\LF{{\rm LF}\,}

\def\P{\mathbb{P}}
\def\t{\mathrm{t}}
\def\i{\mathrm{i}}

\def\s{\mathbf{s}}
\def\a{\mathbf{a}}

\newcommand {\debeq}{\begin{eqnarray*}}
\newcommand {\fineq}{\end{eqnarray*}}

\newtheorem	{thm}		{Theorem}[section]

\newtheorem	{lem} 	[thm]	{Lemma}
\newtheorem     {rem}      [thm]{Remark}
\newtheorem	{prop}	[thm]{Proposition}
\newtheorem	{cor}		[thm]{Corollary}

\begin{document}

\title{The coalescent point process of multi-type branching trees 
}
\author{\textsc{Lea Popovic and Mariolys Rivas}
}
\date{\today}
\maketitle
\begin{abstract}
\noindent

We define a multi-type coalescent point process of a general branching process with finitely many types. This multi-type coalescent fully describes the genealogy of the (quasi-stationary) standing population, providing types along ancestral lineages of all individuals in the standing population. We show that this point-process is a functional of a certain Markov chain defined by the planar embedding of the branching process. 

This Markov chain is further used to determine statistical properties of the ancestral tree, such as the time to the most recent common ancestor (MRCA) of two consecutive individuals in the standing population, as well as of two individuals of the same type. 
These formulae are particularly simple for branching processes with a multi-type \mbox{linear-fractional (LF)} offspring distribution. 
We illustrate their use by showing how an (a)symmetrical offspring distribution affects features of the ancestral tree in a two-type LF branching processes. \end{abstract}  	
\medskip

\noindent \textit{Running head.} The multi-type coalescent point process.\\
\textit{MSC Subject Classification (2000).} \\%Primary 60J80; secondary  60G55, 60G57, 60J10, 60J85, 60J27, 92D10, 92D25.\\
%({\it submit to "Stochastic Processes and their Applications"}??)
\medskip

\tableofcontents

\section{Introduction}

\quad Evolutionary biologists use phylogenetic trees of extant species to study the speciation and extinction patterns of different groups of species.  Numerous phylogenetic trees have been made using sequence data, allowing systematists to investigate the mechanisms of evolutionary dynamics that led to the formation of such a species clade. Most inferences were based on methods that used type independent evolutionary rates of speciation and extinction, and studies of phenotype evolution have not generally accounted for type dependent diversification rates. Recently a number of more type-dependent extensions are being considered (Maddison et al. \cite{MMO07}, Fitzjohnet al. \cite{FMO09}, Igic \& Goldberg \cite{GI08}, \cite{GI12}, etc) with interesting  consequences for validity of previously held beliefs (e.g. Dollo's law).

In the study of macroevolution one uses phylogenetic trees (with or without branch lengths) obtained by genetic sequencing to fit models of speciation and extinction and possibly infer some features of this process, e.g. whether rates are time dependent or trait dependent etc. All (but a small subset) of these methods rely on simulation of branching processes suggested by the model in forward time and calculating the likelihood of obtaining the chosen phylogenetic trees. 
In population genetics an alternative way of assessing the fit of the ancestral relationship based on the phylogenetic trees uses backward simulations based on a coalescent approach. In macroevolution this poses the problem of pre-specifying the unknown random process of fluctuation species numbers over time. One way of avoiding this problem use the point process approach based on a standing branching population as developed in \cite{P04}, further extended by \cite{L10} .

From a mathematical point of view a phylogenetic tree is a genealogical  (or reduced) tree derived from a branching process. The randomness of the branching process implies a distribution on trees that are derived in this manner, and depends on the specifications of the branching process (offspring distributions in discrete time, plus lifetime rates in continuous time). 
A fair amount is known about the distribution of genealogical trees derived from single type branching processes, both for trees of finite size and for the asymptotics as the size becomes infinite.  Much less is known about the same objects for multi-type branching processes, other than asymptotic results based on the almost sure convergence of types for the supercritical multi-type processes (see \cite{KLPP97} and \cite{GB03}). %In \cite{GB03} give almost sure convergence of time and population averages of ancestral types.
 
Deriving an exact distribution for the ancestral tree of a standing population of an arbitrary branching process is unsurprisingly challenging. The first problem, addressed in \cite{LP13}, is handling the branch points in the ancestral tree with multiple surviving offspring. The second problem, arising only for multi-type branching processes, is handling the fact that branch points require the knowledge of the parental type. The latter consequently requires a modification of the contour process approach which keeps track of full memory of the ancestral types of any individual in the standing population (for example via snake construction of \cite{LG93} and \cite{DS00}. We make the analogous modification in the coalescent point process of the ancestral tree.

The unfortunate consequence of this modification is that for the general multi-type branching process some calculations  require an assumption that the branching process is in its the quasi-stationary equilibrium. From the biological point of view this was indeed the most common assumption when reconstructing phylogenies from sequence data, though the logic of making this assumption has recently been questioned in the practice of fitting models of evolution of phenotypes (\cite{GI08}, \cite{FMO09}). From the mathematical point of view, more often than not, quasi-stationary distributions of branching processes can not be calculated explicitly (\cite{M12}), and one would have to resort to numerical calculations for its generating function based on a certain partial difference equation. 

In order to make comparisons on the effect of different branching distributions on the shape of the ancestral trees we get around the problem of quasi-stationarity by focusing on a special case of branching processes. These are multi-type processes whose generating function has a linear-fractional form, and for which many of the point process calculations significantly simplify. A different and complementary set of results on linear fractional multi-type branching processes, on the recurrence on the type space and long term behaviour of such a branching process with countably (rather than finitely) many types, was recently obtained by \cite{Sag13}. A number of their results rely on calculations from the contour process of the multi-type branching process, which is intimately related to our point process construction as well.

\section{Multi-type Coalescent Model and General Results}

\quad In this paper we extend the coalescent point process construction of \cite{LP13} to the case of multi-type branching processes. Our goal is to exploit the Markovian features of the coalescent point process in order to derive features of multi-type phylogenetic trees, and identify the statistics in multi-type phylogenetic trees that are not present in single type trees. We first derive the distribution of the most recent common ancestor of two species from the standing population. We then derive the time of the most recent common ancestor of two species of the same type, and its dependence on the species type. 

\subsection{Multi-type branching process}
\quad We start with notation for multi-type branching processes. 
Let $\{1,2,\ldots,k\}$ be a set of \emph{types} of a population. A \emph{multi-type} or \emph{$k$-type branching process} is a vector-valued Markov process  in discrete time $(\bm{Z}^{(n)})_{n\geq0}$, with $\bm{Z}^{(n)}=(Z_1^{(n)}, Z_2^{(n)},\ldots, Z_k^{(n)})$ a $k$-dimensional random vector whose $\ell$-th coordinate is the number of individuals  of type $\ell$ at generation $n$. Generations will be indexed by $n\in\mathbb{N}_0$ in the superscript, and types will be indexed by $\{1,\ldots,k\}$ in the subscript.

For any $\bm{z}=(z_1,\ldots,z_k)\in \mathbb{N}_0^k$ the matrix  of transition probabilities, and the $n$-th iteration of this matrix,  are denoted by
$$P_{\ell}(\bm{z})=\P(\bm{Z}^{(1)}=\bm{z}\,|\, \bm{Z}^{(0)}=\bm{e}_\ell), \quad P^{(n)}_{\ell}(\bm{z})=\P(\bm{Z}^{(n)}=\bm{z}\,|\, \bm{Z}^{(0)}=\bm{e}_\ell),$$
where $\bm{e}_\ell$ is a unit vector of $\ell$-th coordinate. 
For $\bm{s}=(s_1,\ldots,s_k)\in\mathbb{Z}^k$, the probability generating function of the offspring distribution $\bm{\xi}$ is denoted by $\bm f(\bm s):=(f_1,\ldots,f_k)(\bm s)$ where
$$f_\ell(\bm{s})=\E(\bm{s}^{\bm{Z}^{(1)}}\,|\,\bm{Z}^{(0)}=\bm{e}_\ell)=\sum_{(z_1,\ldots,z_k)\in\mathbb{N}^k} P_{\ell}(\bm{z})s_1^{z_1}\cdots s_k^{z_k},\,\,\,\, \mbox{ for }|s_1|,\ldots,|s_k|\leq1.$$
and the probability generating function of the $n$-th generation population, the $n$-fold composition of $\bm f(\bm s)$, is denoted by $\bm f^{(n)}(\bm s)$ where $f_\ell^{(n)}(\bm{s})=\E(\bm{s}^{\bm{Z}^{(n)}}\,|\,\bm{Z}^{(0)}=\bm{e}_\ell)$. For $n=0$ let $f^{(0)}(\bm s)=\bm s$, and note that $f^{(1)}(\bm s)=f(\bm s)$.

We let 
%$\bm q:=(q_{1},\ldots,q_{k})$ denote the extinction probability, where $\bm 0=(0,\ldots,0)$:
%$$q_{\ell}=P(\exists n: \bm Z^{(n)}=\bm 0\,|\, \bm Z^{(0)}=\bm e_\ell),\,\,\,\,\mbox{for } j=1,\ldots,k$$  and 
$ \bm M=(m_{\ell \ell'})_{1\leq \ell,\ell'\leq k}$ be the matrix of the expected number of offspring of each type from parents of different types:
$$m_{\ell\ell'}=\E(\bm Z^{(1)}_{\ell'}\,|\,\bm Z^{(0)}=\bm e_\ell)=\left.\frac{\partial f_\ell(\bm s)}{\partial s_{\ell'}}\right|_{\bm s=\bm 1},\,\,\,\, \mbox{for } \ell,\ell'=1,\ldots,k,$$
\noindent where $\bm 1=(1,1,\ldots,1)$ and we assume all $m_{\ell\ell'}<\infty$. 
A multi-type Galton-Watson process is called \emph{positive regular} (\emph {or irreducible}) if  for some $n>0$  the mean matrix of its $n$-the generation population $\bm M^n$ is positive (all of the entries $m^{(n)}_{\ell\ell'}>0$ are strictly positive entries).
A process is called \emph{singular} if each individual has exactly one offspring. We assume that the multi-type G-W process is non-singular and irreducible throughout this paper.

As a positive regular matrix $\bm M$ has a maximal, positive and simple eigenvalue $\rho$, which has an associated positive right eigenvector $\bm u$ and a positive left eigenvector $\bm v$. We assume that  $\bm u\cdot \bm v=1$ and $\bm u\cdot \bm 1=1$. The role that $\rho$ plays in the multi-type setting is similar to the role of $\mu$ in the one type case, distinguishing \emph{subcritical, critical or supercritical} processes when $\rho<1$, $\rho=1$ or $\rho>1$, respectively.
\medskip

\subsection{Single-type coalescent point process}
\quad The coalescent point process of a branching tree is a process describing the genealogy of the standing population backwards in time, directly displaying the coalescence times as a sequence running over the
current population size. It constructs a set of points each corresponding to a most recent common ancestors of two individuals in the current population, and its depth (or vertical height) corresponds to the time when the lineages of these two individuals branched off (separated) from each other. The coalescent point process has a bijective correspondence with the ancestral tree of the current population, and allows the full ancestral tree to be reconstructed from its values. It was introduced in \cite{P04} for the ancestral tree of a continuous time single type branching process conditioned on its current population size, and called the {\it genealogical point process}. Its distribution was obtained from its relationship with the {contour} (or {height}) {process} associated with a unit speed traversal of the branching tree. The convenient property of  this particular branching model is that its contour process is Markovian, which implied that the points in this point process are {\it simple} -that is, each branch point had degree two-  and that they are independent samples from the same distribution of depths. This allows one to reconstruct the ancestral tree of a population of $n$ current individuals based simply on a sample of size $n$ from this distribution (see Fig.2 of \cite{P04}). This genealogical point process was used  in \cite{AP05} to obtain statistical information for the ancestral trees of a critical branching process, was extended to non-critical binary processes in \cite{G08} and to homogeneous binary Crump-Mode-Jagers processes in \cite{L10}.

%\begin{figure}\begin{center}
%\includegraphics[totalheight=0.25 \textheight,width=0.7\textwidth]{Fig2[16].pdf}
%\caption{Ancestral tree of $5$ standing individuals (left); its genealogical point-process (right).}
%\end{center}\end{figure}
%\hspace{-5mm}{\it include a representative figure here - genealogical point-process and ancestral tree - Fig2 of \cite{P04}}

This initial construction of the genealogical point process had to be extended to accommodate Galton-Watson branching processes with general  offspring distribution when the contour process of the branching tree is no longer Markovian. In this case depths of points in the process were no longer sufficient in order to fully reconstruct the ancestral tree, as the most recent common ancestors were no longer distinct for every pair of current individuals. In other words, branch points in the ancestral tree no longer always had degree exactly equal to two, and it was necessary to keep track of the multiplicity of these points as given by their branching degree. In \cite{LP13} a construction was made which, rather than having all simple points with mass one, has points with (poisitive) integer valued masses. Each point  again corresponds to a most recent common ancestor of two individuals in the current population, and its depth records the time when the two individuals' lineages separated. The additional mass coordinate of this point records the number of current individuals with the same most recent common ancestor as these two which are embedded after (or horizontally to the right) of them (see Figure 1, taken from \cite{LP13}). This process was called the {\it coalescent point process (with multiplicities)}.

\begin{figure}\begin{center}
\includegraphics[totalheight=0.415 \textheight,width=0.92\textwidth]{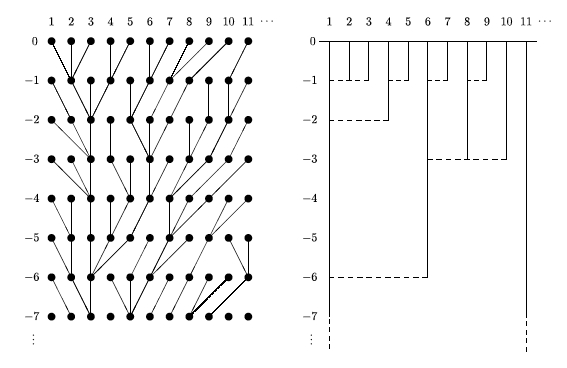}
\caption{Planar embedding of a quasi-stationary branching process (left); coalescence times of consecutive individuals in the standing population (right). The coalescent point process with multiplicities is: $2\cdot\delta_{A_1=1}, 1\cdot\delta_{A_2=1}, 1\cdot\delta_{A_3=2}, 1\cdot\delta_{A_4=1}, 1\cdot\delta_{A_5=6}, 1\cdot\delta_{A_6=1}, 2\!\cdot\delta_{A_7=3}, 1\cdot\delta_{A_8=1}, 1\cdot\delta_{A_9=3}, ...$}
\end{center}\end{figure}
%{\it definitely include another representative figure here - Fig 1 of \cite{LP13}, or one from my slides}

Before we present our extension of this construction we first recall the notation from \cite{LP13}. Consider an arbitrarily large population at the present time from a general quasi-statinory branching process originating at an unspecified arbitrarily large time in the past. In the planar embedding of this process, individuals are located at points of a discrete lattice ($n\in\mathbb{Z}, i\in \mathbb{N}$), where the first coordinate $n$ denotes the generation and the second coordinate $i$ denotes the position of the individual in the planar embedding from left to right. The number of offsprings of individual $(n,i)$ is denoted by $\xi(n,i)$.
The \emph{standing population} is the population at the present time (generation $n=0$), and its \emph{ancestral tree} is the subtree of the branching tree obtained by following only the branches that lead to an individual present in the standing population. The ancestry of an individual  from generation 0 can be traced backwards in time as follows. Define
$$a_i(n):=\mbox{index of the ancestor of individual }(0,i)\mbox{ in generation } -n.$$
 The \emph{coalescent time}  $C_{i,j}$ of individuals $(0,i)$ and $(0,j)$ is the time of the most recent common ancestor between these two, that is, 
 $$C_{i,j}:= \min\{n\geq1\,:\, a_i(n)=a_j(n)\}, \mbox{ with } \min (\emptyset)=\infty.$$
 In particular, define $A_i:=C_{i,i+1}$ {which identifies the coalescent time of individuals $(0,i)$ and $(0,i+1)$.}  It can be easily shown that $C_{i,j}=\max\{A_i,\,A_{i+1},\ldots,A_{j-1}\}$. The sequence $(A_i)_{i\geq1}$ is called the \emph{coalescent point process}. The genealogy back in time of the present population, that is its ancestral tree, is then uniquely determined by the process $(A_i)_{i\geq1}$. (This was sufficient information for the genealogical point process of  binary branching processes in \cite{P04,AP05,G08,L10}.) Define an auxiliary process  $(D_i)_{i\geq1}$ of integer valued sequences $D_i=\{D_i(n), n\geq 1\}$ for each $i\ge 1$, {which records future branch degrees along the ancestral lineage of individual $(0,i)$}
 \begin{eqnarray*}D_i(n)\!\!\!\!&:=\mbox{number of {\it surviving} offspring of individual } (n,a_i(n))\mbox{ in generation } -n\\ & \mbox{ embedded in the {\it ancestral tree} to the right of the lineage of } (0,i) \mbox{ itself}\end{eqnarray*}
It turns out that the process $(D_i)_{i\geq1}$ has all the nice properties needed to identify the law of the coalescent point-process (Theorem 2.1 of \cite{LP13}): $A_i$ is a functional of $D_i$ given by 
$$A_i=\min\{n\ge 1:D_i(n)\neq 0\}$$
and the law of the process $(D_i,i\ge 1)$ is determined by the fact that it is a sequence-valued Markov chain, started at the null sequence $D_0=(0,0,\ldots)$, with transitions given as follows - for any sequence $(d_n; n\ge 0)\in\mathbb{N}^{\mathbb{N}}$ 
$$(D_{i+1}(n)\,|\, D_i(\cdot)=d_{\cdot})\;\mathop{=}^d\;\left\{\begin{array}{lll} d_n&\mbox{ for } & n>A_i,\\
d_{A_i}-1& \mbox{ for } & n=A_i,\\
 \zeta'_n& \mbox{ for  } & 1\leq n<A_i,\end{array}\right.$$
where the random variables $ \zeta'_1, \zeta'_2,\ldots, \zeta'_{A_i-1}$ are independent random variables.

The distributions of variables $\{\zeta'_n\}_ {n\ge 1}$ are specified as follows. 
If $\xi$ is the offspring distribution of this Galton-Watson branching process with probability generating function $f(s)$, the random variables $\xi(n,j)$, representing the number of offspring of individual $(n,j)$ for any indices $n,j\in\mathbb{N}$, are all independent identically distributed as $\xi$. The survival probability to generation $0$ of each offspring of an individual in generation $-n$ is given by $p_{n-1}:=1-f^{(n-1)}(0)$ where $f^{(n-1)}$ is the $(n-1)$-fold composition of $f$. This, in particular, holds for the offspring of $(n,a_i(n))$, the ancestor of $(0,i)$ in generation $-n$.
If we let $\{\epsilon_n^1, \epsilon_n^2, \ldots\}$ be an independent sequence of i.i.d. Bernoulli variables  with parameter $\P(\epsilon_n^m =1)=p_{n-1}$ (we deviate slightly notation from \cite{LP13} here), and use an independent variable $\xi$, we  can define the random sum
$$ \zeta_n: =\sum_{m=1}^{\xi} \epsilon_n^m $$
and, for each $n\ge 1$, the law of $\zeta'_n$ is defined by 
$$ \zeta'_n: \mathop{=}^{d}(\zeta_n-1 \big| \zeta_n\neq 0).$$

\subsection{Multi-type coalescent point process}

\quad Our construction of the coalescent point process for a multi-type Galton-Watson branching tree will rely heavily on the coalescent point process (with multiplicities) and its auxiliary process $(D_i)_{i\ge 0}$. Unlike in the single type branching process where the offspring distribution $\xi$ is the same for all individuals, in order to determine the offspring distribution $\xi(n,i)$ in the multi-type branching process one also has to know the type of the individual $(n,i)$. As a consequence, a Markov process from which the coalescent point process can be reconstructed, will have to contain the information on the individuals' types as well. This, unfortunately, also makes notation for the multi-type process lengthier. Throughout the paper we will reserve boldface symbols for vectors and matrices. 

We start by extending the notation for the planar embedding of a branching process. In the single type case all offspring of an individual were interchangeable in terms of the law of their future subtrees. In the multi-type case we need to make an assumption about  the {\it order of embedding} an individual's offspring of different types in the plane. If there is no particular reason to differentiate the order of the offspring we assume that the order in which they are embedded is chosen {\it uniformly} at random from all possible ways to order them. In the next section we will assume a more {\it specific ordering} in the case where the offspring distribution is linear-fractional.

In addition to its location coordinates each individual has a type associated with it. Let 
$$\t(n,i):= \mbox{{\it type} of the individual }(n,i)$$
Recall that coalescence times between individuals $(0,i)$ and $(0,i+1)$ in generation $0$ are defined as $A_i:= \min\{n\geq1\,:\, a_i(n)=a_{i+1}(n)\}$ for $i\ge 1$, and by convention $A_0=+\infty$. In the multi-type case the ancestral tree also consists of a sequence of types along its lineages: the {{\it ancestral lineage} of individual $(0,i+1)$} back to its most recent common ancestor with individual $(0, i)$ is, for $i\ge 1$, defined as $\bm A_i\in\{1,2\ldots,k\}^{\mathbb{N}_0}$, which includes a special $0$-th coordinate:
$$\bm A_i:=(\t(0,a_{i+1}(0)), \t(-1,a_{i+1}(1)), \ldots, \t(-A_i+1,a_{i+1}(A_i-1))).$$

For a vector $\bm v\in\{1,2\ldots,k\}^{\mathbb{N}_0}$ let $\bm v_{[j]}$ denote its $j$-th coordinate and $\|\bm v\|$ denote its length, with the convention that $ \|\bm v\|=0$ if $\bm v=\emptyset$.
Note that $\|\bm A_i\|=A_i$. Since $a_{i+1}(0)=i+1$, the $0$-th coordinate $\bm A_{i[0]}$ of the vector $\bm A_i$ is the type of the individual $(0,i+1)$. Also, since $A_0=\infty$ the first ancestral lineage $\bm A_0$ consists of types of all individuals on the left most infinite (back into the past) spine of the ancestral tree.

\begin{figure}\begin{center}
\includegraphics[totalheight=0.415 \textheight,width=0.451\textwidth]{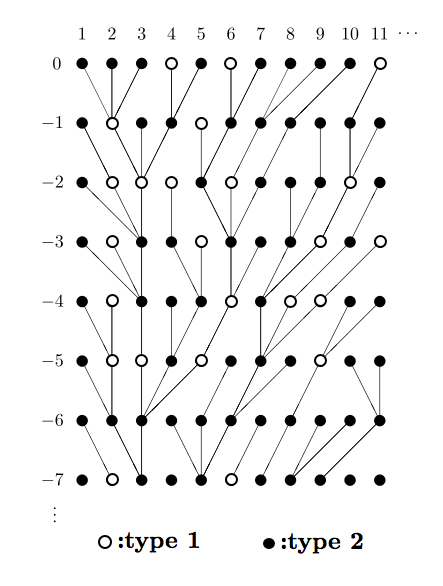}\includegraphics[totalheight=0.415 \textheight,width=0.451\textwidth]{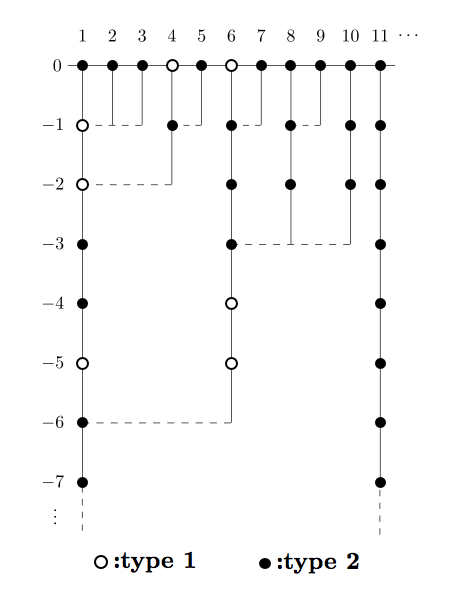}
\caption{Planar embedding of a two-type branching process (left), and the types along the ancestral lineages $\bm A_i$ (right)}
\end{center}\end{figure}

The auxiliary process also needs to be enriched to contain information on types. We define the process $(\bm D_i)_{i\ge 1}$ of vector valued sequences $\bm D_i=\{\bm D_i(n), n\ge 1\}$ in such a way that each $\bm D_i(n)\in\{1,2\ldots,k\}^{\mathbb{N}}$ is a vector of types of offspring of the ancestor $a_i(n)$ in generation $-n$ embedded to the right of the lineage of $(0,i)$ that are `survivors' (meaning that they have progeny that are alive in generation $0$):
 \begin{eqnarray*}\bm D_i(n)\!\!\!\!&:=\mbox{vector of {\it types} of {\it surviving} offspring of individual } (-n,a_i(n))\mbox{ in generation }  -n\\ &\mbox{ embedded in the {\it ancestral tree} to the right of and {\it including} the lineage of } (0,i) \end{eqnarray*}
Clearly $\|\bm D_i(n)\|\ge 1$, and note that $D_i(n):= \|\bm D_i(n)\|-1$ is the number of surviving offspring of individual $(-n,a_i(n))$ embedded to the right of (and excluding) the lineage of $(0,i)$, as in the single type process. 

%As in the single type process, the ancestry of the first individual $(0,1)$, defining the sequence $\bm D_1$, is  the left most spine (infinite back in time) of the genealogical tree of the standing population.

In order to describe the law of $\bm D$ we need to provide notation for surviving lineages. Let $\bm \xi_\ell$ be the offspring distribution of an individual of type $\ell$ with probability generating function $f_\ell(\bm s)$. 
For individual $(n,i)$ the law of the number of its offspring, given that its type is $\t(n,i)=\ell$, is that of $\bm \xi_\ell$. 
The survival probability of an offspring of some type $\ell'$ in some generation $-n'$ is given by $p_{n'-1,\ell'}:=1-f_{\ell'}^{(n'-1)}(0,0,\ldots,0)$ where $f_{\ell'}^{(n'-1)}$ is the $(n'-1)$-fold composition of $f$. 

We consider all the survivor progeny of a generation $-n$ ancestor of some individual from the standing population, and suppose that the type of this generation $-n$ ancestor is $\ell$. 
%{(Later in the construction of $\bm D_{i+1}$ given the value of $\bm D_i$ we will take $n=A_i-1$ and $\ell=\bm D_i(A_i)_{[2]}$; the vector $\bm \eta$ we need to construct will be $A_i-1$ long sequence of vectors in type space initiated by individual of type $\bm D_i(A_i)_{[2]}$)}.
For different $\ell'\in\{1,\ldots, k\}$ let $\{\epsilon^1_{n,\ell'}, \epsilon^2_{n,\ell'}, \ldots\}$ be independent sequences of i.i.d. Bernoulli variables  with parameters $\P(\epsilon^m_{n,\ell'} =1)=p_{n-1,\ell'}$. Start with an independent variable $\bm \xi_\ell$, which takes values in $\mathbb{N}_0^k$ and has $\bm \xi_{\ell, \ell'}$ offspring of type $\ell'$, and define the vector of random sums:
$$\bm \zeta_{n,\ell}: =\big(\sum_{m=1}^{\bm \xi_{\ell, 1}} \epsilon^m_{n,1},\ldots, \sum_{m=1}^{\bm \xi_{\ell, k}} \epsilon^m_{n,k} \big) $$
whose $\ell'$ coordinate is denoted by $\bm \zeta_{n,\ell, \ell'}$.
Then, the law of $\bm \zeta'_{n, \ell}$, which represents the number of surviving offspring  of different types in generation $-(n-1)$ of the initiating generation $-n$ ancestor, is given by: 
$$\bm \zeta'_{n,\ell}: \mathop{=}^{d}(\bm \zeta_{n,\ell} \big| \sum_{\ell'=1}^k\bm \zeta_{n,\ell,\ell'}\neq 0).$$
Let $\bm d(\bm \zeta'_{n,\ell})\in\{1,2,\ldots,k\}^{\mathbb{N}}$ be an ordering of all the offspring counted by $\bm \zeta'_{n,\ell}$ chosen uniformly at random from all possible orderings (or in some specific way, as in the next section). 

Recall that $\bm v_{[j]}$ denotes the $j$-th coordinate of a vector $\bm v\in\{1,2\ldots,k\}^{\mathbb{N}}$. Then $\jmath:=\bm d(\bm \zeta'_{n,\ell})_{[1]}$ is the type of the left most surviving offspring  in generation $-(n-1)$ of the ancestor from generation $-n$. 
Again, for different $\ell'\in\{1,\ldots, k\}$ let $\{\epsilon^1_{n-1,\ell'}, \epsilon^2_{n-1,\ell'}, \ldots\}$ be independent sequences of i.i.d. Bernoulli variables  with parameters $\P(\epsilon^m_{n-1,\ell'} =1)=p_{n-2,\ell'}$ (independent of all earlier sequences of Bernoulli variables).
Proceed with an independent variable $\bm \xi_\jmath$, and define the vector of random sums: 
$$\bm \zeta_{n-1,\jmath}: =\big(\sum_{m=1}^{\bm \xi_{\jmath, 1}} \epsilon^m_{n-1,1},\ldots, \sum_{m=1}^{\bm \xi_{\jmath, k}} \epsilon^m_{n-1,k} \big), \quad \bm \zeta'_{n-1,\jmath}: \mathop{=}^{d}(\bm \zeta_{n-1,\jmath} \big| \sum_{\ell'=1}^k\bm \zeta_{n-1,\jmath,\ell'}\neq 0).$$
and let $\bm d(\bm \zeta'_{n-1,\jmath})\in\{1,2,\ldots,k\}^{\mathbb{N}}$ the ordering of these surviving offspring. Then $\kappa:=\bm d(\bm \zeta'_{n-1,\jmath})_{[1]}$ is the type of the left most surviving progeny in generation $-(n-2)$ of the initiating individual from generation $-n$. 

We proceed in this way recursively until generation $-1$ when we obtain the set of offspring $\bm d(\bm \zeta'_{1,\imath})$. 
In order to collect all types of the left most surviving progeny (and their siblings) in different generations $0,-1, \ldots, -(n-2)$, and $-(n-1)$ in one vector, we define an $n$ long sequence of vectors in type space initiated by individual of type $\ell$: 
$$\bm \eta_{n,\ell}:=\big(\bm d(\bm \zeta'_{1,\imath}), \ldots, \bm d(\bm \zeta'_{n-2,\kappa}), \bm d(\bm \zeta'_{n-1,\jmath}), \bm d(\bm \zeta'_{n,\ell})\big)$$
whose coordinates are then the vectors of surviving offspring types in different generations $\bm \eta_{n,\ell}(1)=\bm d(\bm \zeta'_{1,\imath})$, ... , $\bm \eta_{n,\ell}(n-1)=\bm d(\bm \zeta'_{n-1,\jmath})$, $\bm \eta_{n,\ell}(n):=\bm d(\bm \zeta'_{n,\ell})$.

\begin{figure}\begin{center}
\includegraphics[totalheight=0.425 \textheight,width=0.7\textwidth]{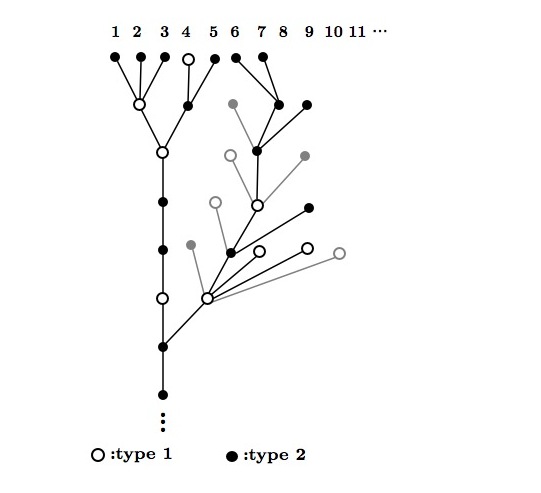}
\caption{Collection of offspring on the ancestral lineage of individual $(0,6)$ originating from its common ancestor in generation $-6$ with individual $(0,5)$ (surviving marked with dark edges;  non-surviving marked with light edges).}
\end{center}\end{figure}
%\hspace{-5mm}{\it include slide 22 from Mariolys' slides but with extra light grey un-surviving offspring}

We introduce one final piece of notation. For a vector $\bm v\in\{1,2\ldots,k\}^{\mathbb{N}}$, let $(\bm v_{[2]},\bm v_{[3]}, \ldots )$ define a vector obtained from $\bm v$ by eliminating the first coordinate $\bm v_{[1]}$ and shifting the rest of its coordinates one coordinate to the left. Having defined the random variables $\bm \eta_{n,\ell}$, for arbitrary $n$ and $\ell$, the reconstruction of the ancestral tree from the auxiliary process is possible as in the single type case.

\begin{figure}\begin{center}
\includegraphics[totalheight=0.35 \textheight,width=1.05\textwidth]{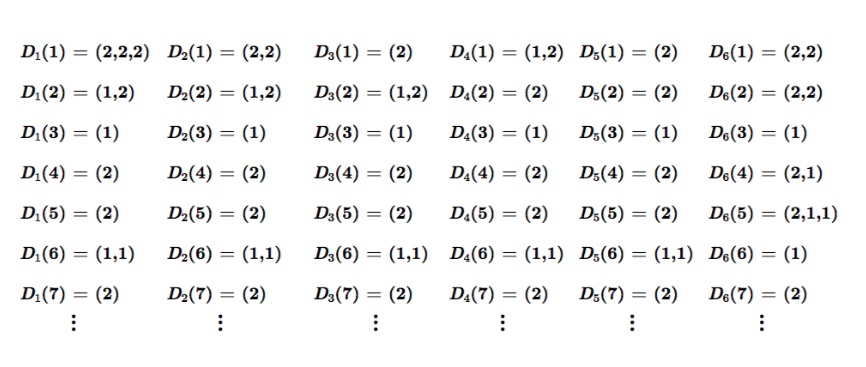}
\caption{Sequences $(D_i(\cdot))_{i\ge 1}$ of surviving offspring types along the lineages of individuals $((0,i))_{i\ge 1}$ corresponding to the ancestral tree given in Figure~2.}
\end{center}\end{figure}
%\hspace{-5mm}{\it include a  figure of $\bm D, i=1..6$ values corresponding the figure above, as in Marilys' Din_1 and Din_2}

\begin{thm}
\label{Din-thm} 
The coalescent times $(A_i)_{i\ge 1}$ and types along the ancestral tree $(\bm A_i)_{i\ge 1}$ are both functionals of $(\bm D_i)_{i\ge 1}$ given by 
$$ A_i=\min\{n\ge 1:\|\bm D_i(n)\|-1\neq 0\}, \quad \bm A_i=\big(\bm{D}_{i+1}(1)_{[1]}, \ldots, \bm{D}_{i+1}(A_i)_{[1]}\big).$$
The sequence $(\bm D_i)_{i\geq1}$ is a Markov chain with transition probabilities given by:
\begin{eqnarray}\label{Din-trlaw}
(\bm{D}_{i+1}(n)\,|\, \bm{D}_i(\cdot)
)\mathop{=}^d\left\{\begin{array}{lll}
\bm{D}_i(n)
&\mbox{\rm for } &n>A_i,\\%({\bm d_n}_{[2]}, {\bm d_n}_{[3]},  \ldots) 
\!\!({\bm{D}}_i(n)_{[2]}, {\bm{D}}_i(n)_{[3]}, \dots)
&\mbox{\rm for }& n=A_i,\\
\bm{\eta}_{_{A_i-1,\bm{D}_{i}(A_i)_{[2]}}}(n) 
&\mbox{\rm for }& 1\le n<A_i \end{array}\right.
\end{eqnarray}
\noindent where the law of the $A_i-1$ long sequence of type vectors $(\bm D_{i+1}(1),\ldots, \bm D_{i+1}(A_i-1))$ is distributed as the vector  $\bm \eta_{n,\ell}$ of types of the left most surviving progeny (and their siblings) in generations $-1,\ldots,-(A_i-1)$  of generation $-n=-(A_i-1)$ individual whose type is $\ell=\bm D_i(A_i)_{[2]}$.
\end{thm}
\textcolor{blue}{\it Proof written \ref{ProofDin-thm} $\Box$}\\

%\textcolor{blue}{\it should we bother defining a process of types at times $(A_i)_{i\ge 0}$, i.e $T_i:=\t(-A_i,a_i(A_i))$?}\\

The Markov chain $(\bm D_{i})_{i\ge 1}$ allows us to calculate some statistical features of the ancestral tree. The most relevant are coalescence times  $(A_{i})_{i\ge 1}$ which indicate the shape of the ancestral tree, and form a non-Markovian process. We next give an  explicit formula for the {\it joint law} of a coalescent time $A_{1}$ of individuals $(0,1)$ and $(0,2)$ {together with the values of types along the whole ancestral lineage $\bm A_0=(\t(0,a_{1}(0)), \t(-1,a_{1}(1)), \ldots)$ of individual $(0,1)$. It illustrates the role of ancestral types when determining branching times in the ancestral tree of the standing population. For a sequence $\bm a\in\{1,\dots,k\}^{\mathbb N_0}$ let $\bm a_{|n'}$ denote the vector of the first coordinates up to $n'$-th one in this sequence 
$\bm a_{|n'}:=(\bm a_{[0]}, \bm a_{[1]}, \dots, \bm a_{[n']})$.

\begin{prop}\label{DistribAi}
For a sequence of types $\bm a=(\bm a_{[0]},\bm a_{[1]},\ldots)\in \{1,\ldots, k\}^{\mathbb N_0}$ 
\begin{eqnarray}\label{A1law}
\P(A_{1}>n,\bm A_{0|n-1}=\bm a_{|n-1} |\,\bm A_{0[n]}=\bm a_{[n]})
=\frac{1}{p_{n,\bm a_{[n]}}}\prod_{n'=1}^n \Bigg(
\frac{\partial f_{\bm a_{[n']}}(\bm s)}{\partial s_{\bm a_{[n'-1]}}} \Big|_{\bm s=\bm 1-\bm p_{n'-1}}\Bigg)
\end{eqnarray}
\noindent where  $\bm 1 -\bm p_{n'-1}:=(1-p_{n'-1,1},\ldots,1-p_{n'-1,k})=\bm f^{(n'-1)}(\bm 0)$ is the vector of extinction probabilities by generation $n'-1$.
\end{prop}
\textcolor{blue}{\it Proof written \ref{ProofDistribAi}\;$\Box$}\\

\noindent An easy modification of the formula above gives $\P(A_{1}=n,\bm A_{0|n-1}=\bm a_{|n-1} |\,\bm A_{0[n]}=\bm a_{[n]})$. We also have the following result on the law of the coalescent time $A_{1}$.
\begin{cor}\label{DistribAiuncond} 
For a single type value $\bm a_{[n]}\in \{1,\ldots, k\}$ and
%\begin{eqnarray}\label{A1law}
%\P(A_{1}>n |\,\bm A_{0[n]}=\bm a_{[n]})
%=\sum_{\bm a_{[\cdot]}\in\bm A_{0[\cdot]}}
%\;\frac{1}{p_{n',\bm a_{[n']}}}\prod_{n'=1}^n \Big(
%\frac{\partial f_{\bm a_{[n']}}(\bm s)}{\partial s_{a_{[n'-1]}}} \Big|_{\bm s=\bm 1-\bm p_{n'-1}}\Big)
%\end{eqnarray}
%where the summation goes over all the possible values of $\bm A_{0[\cdot]}$ for the types of $0\le n'< n$ generations' ancestors of the individual $(0,1)$ which end in $\bm A_{0[n]}=\bm a_{[n]}$. 
%an alternative expression for formula (\ref{A1law}) in terms of the $n$-th generation of 
a multi-type branching process $\bm Z$ 
\begin{eqnarray}\label{A1lawuncond}
\P(A_{1}>n |\,\bm A_{0[n]}=\bm a_{[n]})=\P\Big( \sum_{\ell=1}^kZ^{(n)}_{\ell}=1\,|\, \bm Z^{(n)}\neq \bm 0, \bm Z^{(0)}=\bm e_{\bm a_{[n]}}\Big ).
\end{eqnarray}
\end{cor}
\textcolor{blue}{\it Proof written \ref{ProofDistribAi}\;$\Box$}

%\noindent{\bf Remark}: 
\begin{rem} 
In terms of applications the joint law of $A_1$ and $\bm A_{0|n-1}$ is more useful for reconstructing ancestral trees. Note that the choice of embedding the offspring of each parent uniformly at random in the tree is reflected in these formulae at all. This is in particular evident in (\ref{A1lawuncond}). Moreover, (\ref{A1lawuncond}) can be obtained from (\ref{A1law}) by summing over all the possible values of $\bm A_{0[n']}$ of types of the $0\le n'< n$ generations' ancestors of the individual $(0,1)$ which start with $\bm A_{0[n]}=\bm a_{[n]}$.
In the next section we consider a different choice of ordering the offspring for a specific offspring distribution, and show consistency of these two expressions.\\
\end{rem}

The statistical feature which indicates the distribution of types in the ancestral tree are coalescence times between individuals in the standing population that are of the same type.
Suppose the type of the first individual in generation $0$ is $\bm A_{0[0]}=\ell$, and define the sequence $\i_{\ell,0}:=0, \;\i_{\ell,1}:=\min\{i'> 0: \bm A_{i'[0]}=\ell\},\; \ldots,\;\i_{\ell,i}:=\min\{i'> \i_{\ell,i-1}: \bm A_{i'[0]}=\ell\}, \ldots$ representing the indices of consecutive individuals of type $\ell$ from the standing population.
Define the sequence of \emph{same-type coalescence times} for individuals of type $\ell$ by
 $$B_{\ell, 1}:=\max\{A_{\i_{\ell,1}},\ldots,A_{\i_{\ell,2}-1}\},\;\ldots\;, B_{\ell, i}:=\max\{A_{\i_{\ell,i}},\ldots,A_{\i_{\ell,i+1}-1}\},\ldots$$
Using the values of types on the left most infinite spine $\bm A_0$ we can determine the joint distribution of $B_{\bm A_{0[0]},1}$ and $\bm A_{0|n-1}$ as follows.

\begin{prop}\label{DistribBi}
For a sequence of types $\bm a=(\bm a_{[0]},\bm a_{[1]},\ldots)\in \{1,\ldots, k\}^{\mathbb N}$ with $\bm a_{[0]}=\ell$
\begin{eqnarray}\label{B1law}
\P(B_{\ell,1}>n, \bm A_{0{|n-1}}=\bm a_{|n-1}|\,\bm A_{0[n]}=\bm a_{[n]}, \bm A_{0[0]}=\ell)
=\frac{1}{p_{(n,\ell),\bm a_{[n-1]}}}
\!\!\prod_{n'=1}^n \Big(\frac{\partial f_{\bm a_{[n']}}(\bm s)}{\partial s_{\bm a_{[n'-1]}}} \Big|_{\bm s=\bm 1-\bm p_{(n'-1,\ell')}}\Big)
\end{eqnarray}
\noindent where  $\bm 1-\bm p_{(n'-1,\ell)}:=(1-p_{(n'-1.\ell),1},\ldots,1-p_{(n'-1,\ell), k})=\bm f^{(n'-1)}(\bm {\hat e}_{\ell'})$  with $\bm {\hat e}_{\ell'}=\bm 1-\bm e_{\ell'}$, is the vector of extinction probabilities for lineages with type $\ell$ descendants after  $n'-1$ generations. 
%The summation goes over all possible values of $\bm A_{0[\cdot]}$ for types of $0\le n'< n$ generations' ancestors of the individual $(0,1)$ with the initial value $\bm A_{0[0]}=\ell$ and end value $\bm A_{0[n]}=\bm a_{[n]}$.
\end{prop}
\textcolor{blue}{\it Proof written \ref{ProofDistribBi}\;$\Box$}\\

\noindent As before, we also have the following result on the law of the same-type coalescent time $B_{\ell,1}$.

\begin{cor}\label{DistribBiuncond} 
For a single type value $\bm a_{[n]}\in \{1,\ldots, k\}$ and a multi-type branching process $\bm Z$
\begin{eqnarray}\label{B1lawuncond}
\P(B_{\ell,1}>n |\,\bm A_{0[n]}=\bm a_{[n]},\bm A_{0[0]}=\ell)=\P\big(Z^{(n)}_{\ell}=1\,|\, Z^{(n)}_{\ell}\neq 0, \bm Z^{(0)}=\bm e_{\bm a_{[n]}}\big ).
\end{eqnarray}
\end{cor}
\textcolor{blue}{\it Proof written \ref{ProofDistribBiuncond}\;$\Box$}

\vspace{5mm}It might be tempting to provide a formula for the distribution of $A_2,A_3,\dots$ and $B_{\ell,2},B_{\ell,3},\dots$ in a similar vein using the values of the types on the ancestral lineage of the individuals $(0,2), (0,3),\ldots$ respectively. The information, analogous to that of types of individuals on the left most (infinite into the past) spine $\bm A_0$ used in the formulae for $A_1$ and $B_{\ell,1}$,
which one would need to use for $A_{i+2}$ would be the types along the ancestral lineage of $(0,i+1)$. That is, for $i\ge 0$ one could define the {\it infinite ancestral lineage} $\bm A^\infty_{i}$  of individual $(0,i+1)$ as the infinite sequence:
$$\bm A^\infty_{i}:=(\t(0,a_{i+1}(0)), \t(-1,a_{i+1}(1)), \t(-2,a_{i+1}(2)),\ldots),$$
Note that the restriction of $\bm A^\infty_{i}$ to its first $A_{i}$ entries equals the sequence $\bm A_{i}$ called the {\it ancestral lineage} of $(0,i+1)$, and that $\bm A^{\infty}_0=\bm A_0$.
It is easy to see, as a consequence of Theorem~\ref{Din-thm},  that $(\bm A^\infty_i)_{i\ge 1}$ is also a functional of $(\bm D_i)_{i\ge 1}$ given by:
$$\bm A^\infty_i=(\bm D_{i+1}(1)_{[1]}, \bm D_{i+1}(2)_{[1]}, \bm D_{i+1}(3)_{[1]},\ldots).$$
This follows from the fact that the first $1\le n\le A_i$ entries in this sequence are the same as in $\bm A_i$, while for the subsequent $n>A_i$ entries the ancestry of $(0,i+1)$ is equal to the ancestry of $(0,i)$ (as their ancestors already coalesced) and we have $\bm D_{i+1}(n)=\bm D_i(n)$.\\

Now, providing a formula for $\,\P(A_{i+2}>n|\,\bm A_{i+1[n]}=\bm a_{[n]})\,$ for any $i\ge 0$ can be done only in case the branching mechanism is such that in the coalescent point-process all points are simple (have multiplicities equal to one). This is because, in case of multiple coalescence points, all ancestral lineages, except for that of $(0,1)$, in addition to information about the lineage of individual $(0,i)$ also contains information about the ancestral lineages of $(0,i'), 1\le\i'<i$. In other words, the calculation (used in the proofs of the Propositions~\ref{DistribAi} and \ref{DistribBi}, see Section~\ref{AllProofs}) which exploits the equivalence $\{A_1>n,\bm A_{0|n-1}=\bm a_{|n-1}\}$ iff  \{individuals $\bm a_{[n']}, 1\le n'\le n$ on the ancestral lineage $\bm A_0$  have a single offspring with descendants surviving to generation $0$\} is no longer valid in general. This is a clear  consequence of the fact that the process $(\bm A_i, i\ge 0)$ is itself not Markovian, unless we are in the special case for the offspring distribution which is of linear-fractional form. In the next Section we explore this special case, and extend the above results for the coalescent times as well as for the same-type coalescence times.

\begin{rem}
The construction of the ancestral tree is based on the Markovian property of the auxiliary process $(\bm D_i)_{i\ge 1}$. To initiate the process, we need to draw $\bm D_1$ from $\bm \eta_{\infty, \ell_\infty}$ for some (infinitely old) originating type $\ell_\infty$. One way to draw from this distribution would be to start from generation $0$ and time reverse the quasi-stationary distribution for the branching process. For the case of a single type branching process a construction like this was discussed in \cite{G99}. The multi-type extension is straightforward, with the infinitely old originating individual having, in the $\rho\le 1$ case, the size-biased version of the offspring distribution  given by
$$\widehat \P(\bm Z^{(1)}=\bm{z}|\bm Z^{(0)}=\bm{e}_\ell)=\frac{\P(\bm Z^{(1)}=\bm{z}|\bm Z^{(0)}=\bm{e}_\ell)}{\rho}\frac{\bm{z}\cdot u}{\bm{e}_\ell\cdot u} $$.
\end{rem}

\section{Special case: Linear-fractional branching processes}

Many of the complications which arise in calculating the distribution of ancestral trees in multi-type branching processes simplify a great deal in the special case when the offspring distribution is of the linear-fractional (LF) type. This type of  offspring distribution leads to a number of particularly nice features involving the memoryless property of the geometric distribution. We first recall the definition of the multi-type linear-fractional offspring distribution, and then give a series of specific results for the distribution of the ancestral tree of the standing population, which both illustrate and extend our general results form the previous Section.

\subsection{Multi-type linear-fractional branching process}

We use the same notation as \cite{Sag13} for ease of drawing on known results and making comparisons. Let $\bm H$ be a $k\times k$ sub-stochastic matrix, that is, each row $\bm h_\ell$ of $\bm H$ is a non-negative vector with $\sum_{\ell'=1}^k h_{\ell\ell'}\le 1$, and let $h_{\ell 0}=1-\sum_{\ell'=1}^k h_{\ell\ell'}$. Let $\bm g$ be a non-negative vector such that $\bm g \bm 1^{\t}=\sum_{\ell'=1}^k g_{\ell'}=1$. Let $m> 0$. For any $\bm z=(z_1,\ldots, z_k)\in\mathbb{N}_0^k$ let $|\bm z|=\sum_{\ell=1}^k z_\ell$.

A random vector $\bm \xi_\ell$ taking values in $\mathbb{N}_0^k$ has a \emph{linear-fractional distribution \LF$(\bm h_\ell,\bm g,m)$} if for any  non-negative integer vector $\bm z=(z_1,\ldots, z_k)$ 
$$\P\big(|\bm \xi_\ell|=0\big)=h_{\ell 0},\quad \P\big(\bm \xi_\ell=\bm e_{\ell'}+\bm z\big)=h_{\ell \ell'} \frac{m^{|\bm z|}}{(1+m)^{|\bm z|+1}}\left(\begin{array}{cc}  |\bm z|\\ z_1,\ldots, z_k\end{array}\right)\bm g^{\bm z},$$
where $\bm g^{\bm z}=g_1^{z_1}\cdots g_k^{z_k}$.
The probability generating function of $\bm \xi_\ell$ has the linear fractional form 
$$f_\ell(\bm s)=h_{\ell 0}+\frac{\sum_{\ell'=1}^kh_{\ell \ell'}s_{\ell'}}{1+m-m\sum_{\ell'=1}^k g_{\ell'}s_{\ell'}}.$$ 
One can also represent the random vector $\bm \xi_\ell$ as a sequence of offsprings, where the first offspring has type distribution given by $\bm h_\ell$ and the children after the first one have geometric distribution with mean $m$ and type distribution given by $\bm g$ independently for each offspring. Moreover, the probability generating function of $(\bm \xi_\ell |\bm \xi_\ell \neq \bm 0)$ is that of a shifted Multivariate-Geometric distribution 
$$\E(\bm s^{\bm \xi_\ell} |\bm \xi_\ell \neq \bm 0)=\frac{(1-h_{\ell 0})^{-1}\sum_{\ell'=1}^kh_{\ell \ell'}s_{\ell'}}{1+m-m\sum_{\ell'=1}^k g_{\ell'}s_{\ell'}}.$$

A {\it multi-type linear-fractional branching process} $\LF\!(\bm H, \bm g, m)$ is a branching process in which each individual of type $\ell\in\{1,\ldots,k\}$ reproduces according to the  $\LF\!(\bm h_\ell,\bm g, m)$ offspring distribution $\bm \xi_\ell$. In other words, the probability generating function of the offspring distribution is $\bm f(\bm s)=(f_1,\ldots,f_k)(\bm s)$ with $f_\ell$ as above. Its mean matrix is given by $\bm M=\bm H+m\bm H\bm 1^{\t}\bm g$.

The fact that the parameters $\bm g$ and $m$ do not depend on the parent's type ensures that the population size in each generation of this process also has a linear fractional distribution. This is given by the following known result.

\begin{thm}[Proposition 3 \cite{JL06}, Theorem 3 \cite{Sag13}]\label{thm-LF}
The $n$-th generation population size vector $\bm Z^{(n)}$ of a multi-type linear-fractional branching process $\LF\!(\bm H,\bm g, m)$ started with one individual $\bm Z^{(0)}=\bm e_\ell$ has a linear-fractional distribution $\LF\!(\bm h_\ell^{(n)},\bm g^{(n)},m^{(n)})$ whose parameters are determined by:
\begin{eqnarray}\label{LFparam}
m^{(n)}\!\!\!&=&\!\!m\bm g (\bm I+\bm M+\cdots +\bm M^{n-1}) \bm 1^{\t},\nonumber\\ 
\bm g^{(n)}\!&=&\!\!\!\!\frac{\!\!\!m}{\,m^{(n)}}\,\bm g(\bm I+\bm M+\cdots +\bm M^{n-1}),\\ 
\bm H^{(n)}\!\!\!&=&\!\!\bm M^n-\frac{m^{(n)}}{1+m^{(n)}}\bm M^n\bm 1^t\bm g^{(n)},\nonumber
\end{eqnarray}
where the vector $\bm h_\ell^{(n)}$ is the $\ell$-th row of the matrix $\bm H^{(n)}$, and $\bm 1^{\t}$ is the transpose of $\bm 1=(1,\ldots,1)$. 
\end{thm}
Note that, as a consequence, $(\bm Z^{(n)} |\bm Z^{(n)} \neq \bm 0, \bm Z^{(0)}=\bm e_{\ell})$ is distributed as a shifted Multivariate-Geometric distribution 
$$\E(\bm s^{Z^{(n)}} |\bm Z^{(n)}\neq \bm 0, \bm Z^{(0)}=\bm e_{\ell})=\frac{(1-h_{\ell 0}^{(n)})^{-1}\sum_{\ell'=1}^kh_{\ell \ell'}^{(n)}s_{\ell'}}{1+m^{(n)}-m^{(n)}\sum_{\ell'=1}^k g_{\ell'}^{(n)}s_{\ell'}}.$$

This theorem was proved in \cite{JL06} using an algebraic approach, while \cite{Sag13} provided a different proof using the jumping contour representation of the branching process and its nice Markovian structure. As shown in \cite{Sag13} this result and the linear-fractional definition also apply for the case of infinitely many types. See \cite{Sag13} also for a discussion about the use of linear-fractional branching processes in applications.

\subsection{Coalescent times in linear-fractional ancestral tree}

We use the coalescent point-process construction to get simpler results for the distribution of coalescent times (and same-type coalescent times) for this special class of multi-type branching processes. For this purpose we make one change in our original construction pertaining to the embedding of the multi-type tree in the plane. For a general offspring distribution we made the assumption that the offspring of any parent are embedded in a left to right order chosen {\it uniformly} at random from all possible orderings. For the linear-fractional offspring distribution we make a particular assumption that the offspring with distribution given by the vector $\bm h_\cdot$ is embedded as the left most individual, followed by the rest of the offspring according to an arbitrary order. 

\begin{prop} \label{DistribAiLFCase2}
The coalescence times $(A_i)_{i\ge 0}$ in the ancestral tree of a $\LF\!(\bm H, \bm g,m)$ branching process are independent identically distributed variables with
\begin{equation}\label{A1lawLFCase2}
P(A_1>n)=\prod_{n'=1}^n \frac{1}{1+m-m\sum_{\ell'=1}^k \bm g_{\ell'} h_{\ell'0}^{(n'-1)}}=\frac{1}{1+m^{(n)}},
\end{equation}
where $h_{\ell'0}^{(n'-1)}\!\!=1-\bm h_{\ell'}^{(n'-1)}\bm 1^{\t}$, $\bm h_{\ell'}^{(n'-1)}\!\!$ is the $\ell'$-th row of the matrix $\bm H^{(n)}$ from (\ref{LFparam}), $h^{(0)}_{\ell'0}=0\; \forall \ell'$.
The law of the coalescent times also satisfies
$\P\big(A_{1}>n |\,\bm A_{0[n]}=\bm a_{[n]}\big)=\P(A_1>n)$.
\end{prop}
\textcolor{blue}{\it Proof written \ref{ProofDistribAiLFCase2}\;$\Box$}

\begin{rem} The second, simpler expression in (\ref{A1lawLFCase2}) can be obtained either using arithmetic properties of parameters (\ref{LFparam}) from Theorem~\ref{thm-LF}, or using the expression in terms of the $n$-th generation of a multi-type LF branching process.
\end{rem}

\begin{prop} \label{DistribBiLFCase2}
For any type $\ell\in\{1,\dots,k\}$, the same-type coalescence times $(B_{\ell,i})_{i\ge 0}$ are independent identically distributed variables with 
\begin{equation}\label{B1lawLFCase2}
P(B_{\ell,1}>n|\,\bm A_{0[0]}=\ell)=\prod_{n'=1}^n \frac{1}{1+m-m\sum_{\ell'=1}^k \bm g_{\ell'} \tilde h_{\ell'0}^{(n'-1)}}=\frac{1}{1+m^{(n)}g^{(n)}_{\ell}}
\end{equation}
where, for $n'>1$, $\tilde{\bm h}_{\ell'0}^{(n'-1)}$ is given in terms of $h_{\ell'0}^{(n'-1)}=1-\bm h_{\ell'}^{(n'-1)}\bm 1^{\t}$ and parameters in (\ref{LFparam}) by
$$ \tilde h_{\ell' 0}^{(n'-1)}=h^{(n'-1)}_{\ell'0}+\frac{1-h^{(n'-1)}_{\ell'0}-h^{(n'-1)}_{\ell'\ell}}{1+m^{(n'-1)}g^{(n'-1)}_{\ell}},  
$$
and $\tilde h_{\ell'0}^{(0)}=1\; \forall \ell'\neq\ell$, while $\tilde h_{\ell 0}^{(0)}=0$. 
\end{prop}
\textcolor{blue}{\it Proof written \ref{ProofDistribBiLFCase2}\;$\Box$}\\

Although the multi-type LF offspring distribution has a seemingly small level of dependence between the offspring and parent type, it still affects the distribution of types in the ancestral tree. One can also consider a multi-type branching process where offspring distribution is completely independent of the parent type. In this case the shape of the tree and the types on the tree can be decoupled, and the distribution of types is only governed by the frequency of this type in the population.

Consider a special case of a LF distribution where $\bm H:=\bm 1^{\t}\bm h$, for $\bm h\bm 1^{\t}=h_1+\cdots+h_k\le 1$,  then each parent has the same $\LF\!(\bm h, \bm g, m)$ offspring distribution. Further, if we are to have no distinction between the first offspring and the rest, then consider $\bm h:=(1-h_0)\bm g$, for $h_0\in(0,1)$. In this case all parents have the same offspring laws where their number of offspring has a single-type LF$(h_0,m)$ distribution with probability generating function $h_0+(1-h_0)s/(1+m-ms)$, given the number of their offspring, the distribution of types for these offspring is Multinomial with parameter $\bm g$. In this case we get the following formulae for the law of coalescence times and same-type coalescence times.

\begin{cor}\label{DistribLFCase1}
If the offspring distribution of each parent is independent of the parent's type with $\LF\!((1-h_0)\bm 1^{\t}\bm g, \bm g, m)$ distribution, then 
\begin{eqnarray*}\label{A1lawLFCase1}
\P(A_1>n)&=& \Bigg\{\begin{tabular}{ll}$\frac{m-h_0(1+m)}{m(1+m)^n(1-h_0)^n-h_0(1+m)}$&$\mathrm{if}$\; $(1-h_0)(1+m)\neq 1$\\ $ $\\
 $\frac{1-h_0}{1-h_0+nh_0}$\;\;\;\qquad\qquad\qquad\qquad\qquad &$\mathrm{if}$\; $(1-h_0)(1+m)=1$\end{tabular}
\end{eqnarray*}
and
\vspace{-2mm}
\begin{eqnarray}\label{B1lawLFCase1}
\P(B_{\ell, 1}>n|\,\bm A_{0[0]}=\ell)=\frac{1-\P(A_1\le n)}{1-\P(A_1\le n)(1- g_\ell)}%=\frac{\P(A_1> n)}{\P(A_1> n)+(1-\P(A_1> n))g_\ell}
\end{eqnarray}
\end{cor}
\textcolor{blue}{\it Proof written - \ref{ProofDistribLFCase1}$\Box$}\\

As expected, the distribution of types $\bm g$ has no effect on the law of $(A_i)_{i\ge 1}$ and the shape of the tree, but appears in the distribution of types in the tree as indicated by the law of  $(B_{\ell,i})_{i\ge 1}$. 

We can consider the process of coalescent times $(A_i)_{i\ge 0}$ as a simple point-process $A$ on $\{1,2,\dots\}\times\{-1, -2, \dots\}$ with intensity measure $\nu_A\big[\{i\}\times\{-(n+1),\dots\}\big]=\P(A_1>n)$, $\forall i\ge 1$. Similarly  for each $\ell\in\{1,\dots,k\}$, $(B_{\ell,i})_{i\ge 0}$ can be regarded as a simple point-process $B_\ell$ with intensity $\nu_{B_\ell}\big[\{i\}\times\{-(n+1),\dots\}\big]=\P(B_{\ell,1}>n)$ on $\{1,2,\dots\}\times\{-1, -2, \dots\}$. Note that for any $i\ge 1$, $n\ge 1$ (\ref{B1lawLFCase1}) implies that
$$
\nu_{B_\ell}\big[\{i\}\times\{-1,\dots,-n\}\big]=\frac{\nu_A\big[\{i\}\times\{-1,\dots,-n\}\big]g_\ell}{\nu_A\big[\{i\}\times\{-(n+1),\dots\}\big]+\nu_A\big[\{i\}\times\{-1,\dots,-n\}\big]g_\ell}
$$
showing that only a fraction of all coalescent times are candidates for same-type coalescence times for type $\ell$. 
Intuitively, when considering $B_{\ell,i}$ mark each coalescent time $A_i,A_{i+1},\dots$ with the probability that the next individual in the standing population is of type $\ell$, which is $g_\ell$.
Then, from the filtered view of $B_{\ell,i}$, a coalescence time $A_i,A_{i+1},\dots$ either occurs outside the set $\{1,\dots,n\}$, or it occurs inside this set and it links to a standing individual of type $\ell$.

Note that the intensity measures $\{\nu_{B_\ell}\}_{\ell\in\{1,\dots,k\}}$ do not partition in full the measure $\nu_{A}$, since for any $n\ge 1$ such that  $\P(A_i>n)>0$ we have that
$$
\sum_{\ell=1}^k\nu_{B_\ell}\big[\{i\}\times\{-1,\dots,-n\}\big] <\sum_{\ell=1}^k {\nu_A\big[\{i\}\times\{-1,\dots,-n\}\big]g_\ell}={\nu_A\big[\{i\}\times\{-1,\dots,-n\}\big]}.
$$
This is a consequence of the fact that not all coalescence times are in fact same-type coalescence times for some $\ell$ (for example, in Figure~2 the coalescence time $A_4=1$ of $(0,4)$ and $(0,5)$ is neither a same-type $\bm 1$ nor a a same-type $\bm 2$ coalescence time).

%This also provides the following relationship between the intensity measures $\{\nu_{B_\ell}\}_{\ell\in\{1,\dots,k\}}$.  For fixed $n\ge 1$ consider the index $\;\i_{\ell,n}=\inf\{i\ge 1: \nu_{B_\ell}\big[\{1,\dots, i\}\times \{-(n+1),\dots\}\big]\neq \emptyset\}$ of the first point of $B_\ell$ that falls outside of $\{-1,\dots, -n\}$. Since each point-process $B_\ell$ is simple the distribution of $\i_{\ell,n}$ is Geometric with parameter $\nu_{B_\ell}\big[\{i\}\times\{-(n+1),\dots\}\big]=\P(B_{\ell,1}>n)$ with expected value  $\E[\i_{\ell,n}]=\P(B_{\ell,1}> n)/\P(B_{\ell,1}\le n)$. 

\subsection{Comparison of ancestral trees in two-type models}

We next give an example of using the same-type coalescent times to investigate the effect of differences in offspring distribution on the distribution of types in the ancestral trees they produce. One question that motivated our work is the effect of different diversification rates for different types of individuals (phenotypes). We translate these questions into a discrete time defining an asymmetrical offspring distribution law. 

Specifically, in a population with only two types of individuals, if the transition rates of one type to the other are relatively high, while the other type never transitions into the first, this will be reflected in the distribution of types along the tree. In a discrete time process this is translated in the probability of a parent of the first type giving birth to individuals of the second type and vice versa. We consider this difference in the context of a two-type LF offspring distribution. In order to investigate only the effect on the distribution of types, we will make the distribution of the shape of the tree the same in both cases.

We consider the following two LF offspring distributions on $k=2$ types of individuals. Let the parameters $\bm g=(g, 1-g), g\in[0,1/2]$, $m>0$ and $h_1=1-h_0, h_0\in [0,1]\;$ be the same in both distributions, and for $p\in(0,1)$ let %\textcolor{blue}{\it notation change: $b\mapsto g, S(\cdot)\mapsto G(\cdot)$}
$$\bm H_{\s}=h_1\left(\begin{array}{cc} p  &1-p\\ 1-p& p\end{array}\right),\quad\bm H_{\a}=h_1\left(\begin{array}{cc} p  &1-p\\ 0& 1\end{array}\right),$$
be,  respectively, associated with the symmetrical and the asymmetrical offspring distribution. 
In the symmetrical case parents of either type produce offspring of their own type and of the other type. In the asymmetrical case only a parent of type $\bm 1$ will do that, while a parent of type $\bm 2$ can only produce offspring of its own type.  
Since the number of offspring of each parent depends only on $h_0$ and $m$, the distribution of the ancestral tree with types erased will be the same in both cases. However, the distribution of the two types $\bm 1$ and $\bm 2$ are different, as can be seen in the following result.

\begin{rem}Note that we can assume without loss of generality that $g\in[0,1/2]$, since in case $g\in[1/2,1]$ we can simply reverse the notation of the two types. For $p=1$ there is no asymmetry, nor are there offspring of different type than the parent - individuals in the whole tree are all of the same type. %For $p=0$ the asymmetry is very pronounced - in the symmetric the first born individuals is always of different type than the parent, while in the asymmetric case the first born individual is always of type 2. 
For $(g,p)=(1/2,1/2)$ the symmetric case is special, and the offspring distribution is independent of the type of the parent, as discussed in  Corollary~\ref{DistribLFCase1}.
\end{rem}

\begin{prop}\label{2typecomp}
The distributions of coalescence times $(A_i)_{i\ge 1}$ are the same in both cases. While the distribution of same-type coalescence times $(B_{\bm 1,i})_{i\ge 1}$ and $(B_{\bm 2,i})_{i\ge 1}$ satisfy the following stochastic dominance relations: $\forall p\in [0,1]$, 
$$\P_{\a}(B_{\bm 1,i}>n\,|\,\bm A_{0[0]}=\bm 1)\geq \P_{\s}(B_{\bm 1,i}>n\,|\,\bm A_{0[0]}=\bm 1),$$ and
$$\P_{\a}(B_{\bm 2,i}>n\,|\,\bm A_{0[0]}=\bm 2) \leq \P_{\s}(B_{\bm 2,i}>n\,|\,\bm A_{0[0]}=\bm 2).$$
Also $\forall p\geq 1/2$  the two above inequalities are related by:
$$\P_{\s}(B_{\bm 1,i}>n\,|\,\bm A_{0[0]}=\bm 1)\geq \P_{\s}(B_{\bm 2,i}>n\,|\, \bm A_{0[0]}=\bm 2).$$
\end{prop}
\textcolor{blue}{\it Proof written - \ref{Proof2typecomp}\;$\Box$}\\

The explicit formulae for all of the above probabilities in terms of the parameters $g, h_1, p$ are complicated and can be found in the proof of the Proposition. We see that the consequence of asymmetry (irrespective of the value of $p$) is that the same-type coalescence times are typically shorter for type $\bm 2$ than in the symmetrical case, while they are longer for type $\bm 1$. This intuitively make sense, since subtrees of a type $\bm 2$ can  only contain type $\bm 2$ individuals, while subtrees of a type $\bm 1$ individual contain a mixture of types. 

We can also see the effect that the `strength' $p$ of not transitioning to a different type plays in the symmetric case. When $p\ge 1/2$ having the same type offspring as parent is more likely. In the symmetric case $g\le 1/2$ further implies that type $\bm 1$ is overall less frequent than type $\bm 2$ in the tree. Hence, one would expect that the same-type coalescence times are typically going to be longer for type $\bm 1$ than for type $\bm 2$. \\
%\textcolor{blue}{\it anything else we want to say using formula for $G$ here??}

\section{Discussion of Results}\label{Disc}

With this work we provided the following: 

(1) an explicit and algorithmic way to construct an ancestral tree of the standing population of a (quasi-stationary) multi-type branching process in terms of a Markov chain; and

(2) explicit formulae for calculating: (2a) the basic statistical features that describe the ancestral tree (the law of coalescence times together with the types on the ancestral lineages), as well as (2b) statistical features that link types in the standing population with the shape of the tree (the law of same-type coalescence times). 

As an example of what one can infer from these results, we considered  the special case of a multi-type branching process with linear-fractional offspring distribution, and we obtained very simple formulae for these two sets of statistical features. 
These formulae were then used to assess the differences in the ancestral trees of two different linear-fractional offspring distributions: one `symmetrical' in its treatment of different offspring types, and the other completely `asymmetrical' in that sense. The `symmetry' and `asymmetry' were clearly featured in the statistics of the ancestral trees, which could be used  to infer the extremeness of parameters that determine this (a)symmetry in the offspring distribution.

\section{Proofs}\label{AllProofs}

\subsection{Proof of results for general multi-type coalescent point process}

We first state a spine decomposition of a multi-type branching process conditioned on survival to a certain generation, which shows that, if we consider the infinite (back into the past) lineage of a current individual, at every generation back in the past the subtrees of siblings of the ancestor in that generation are independent of the infinite lineage and are distributed as trees of an unconditioned multi-type branching process. Moreover, {\it knowing the values of their own initial individuals}, these trees are independent from their sibling subtrees, and are independent of their rank in the planar ordering.

For single-type processes this result first appeared in \cite{LPP95} and \cite{G99}. 
For  multi-type processes a decomposition of a tree relative to a spine that is infinite {\it  into the future} is stated in \cite{KLPP97}, and in \cite{GB03} for branching in continuous time. We present a statement in the form of Lemma 2.1 from \cite{G99} for decomposition of trees conditioned only to survive to a fixed generation, and give its proof. Consider a multi-type branching process $\bm Z$ which is still non-extinct in generation $n+1$, let $T^{(i)}, 1\le i\le  |\bm Z^{(1)}|$ denote the subtrees descending from the offspring in the first generation. Let $\bm d(\bm{Z}^{(1)})$ be a uniform ordering of all the offspring types in the first generation, and let $R_{n+1}$ be the rank of the first offspring whose descendants survive to generation $n+1$. 

\begin{lem}\label{GeigMultitype}
The subtrees $T^{(i)}$, $1\leq i\leq  |\bm Z^{(1)}|$, $d_i\in\{1,2,\ldots,\,k\}$, are conditionally independent given $\{\bm{Z}^{(0)}=\bm e_\ell, \bm{Z}^{(1)}=\bm{z}, \bm d(\bm{Z}^{(1)})=\bm d(\bm z),{R}_{n+1}=j \}$, for $1\leq j\leq |z|$ and $\bm{z}=(z_1,\,z_2,\ldots,\,z_{k})$ with $\bm d(\bm z)=(d_1,\ldots, d_{|\bm z|})$ 
\begin{eqnarray*}
\big(T^{(i)}\, | \,  \bm{Z}^{(0)}=\bm e_\ell,\bm{Z}^{(1)}=\bm{z},\bm d(\bm z)\!\!\!\!\!\!&=&\!\!\!\!\!\! (d_1,\ldots, d_{|\bm z|}),{R}_{n+1}=j \big)\\
&&\mathop{=}^d\left\{\begin{array}{l} \big(\mathcal{T}\big|\, \bm Z^{(n)}(\mathcal T)=\bm 0, \bm{Z}^{(0)}(\mathcal T)=\bm e_{d_i}\big),\, 1\leq i\leq j-1\\
 \,\\
\big(\mathcal{T}\big|\, \bm Z^{(n)}(\mathcal T)\neq \bm 0, \bm{Z}^{(0)}(\mathcal T)=\bm e_{d_i}\big),\, i=j\\
\,\\
\big(\mathcal{T}\big|\, \bm{Z}^{(0)}(\mathcal T)=\bm e_{d_i}\big),\, j+1\leq i\leq |\bm z|,\end{array}\right.
\end{eqnarray*}

\noindent where $\mathcal{T}$ denotes the law of a tree of multi-type type branching processes with the p.g.f. of $\bm Z$. 
Further, the conditional joint distribution of ${R}_{n+1}, \bm{Z}^{(1)},\bm d(\bm Z^{(1)})$ is given by
\begin{eqnarray*}\label{distriRZ}
&&\!\!\!\!\!\!\!\!\!\P({R}_{n+1}=j, \bm{Z}^{(1)}=\bm{z}, \bm d(\bm{Z}^{(1)})=\bm d(\bm z)\,\big|\, \bm Z^{(n+1)}\neq \bm 0, \bm{Z}^{(0)}=\bm e_\ell)\\ &&= \frac{{\P(\bm \xi_\ell= \bm z)\P(\bm d(\bm z)=(d_1,\ldots, d_{|\bm z|}))\P(\bm Z^{(n)}\neq \bm 0| \bm{Z}^{(0)}=\bm e_{d_j})\,}\prod_{\ell'=1}^{j-1}\P(\bm{Z}^{(n)}=\bm{0})|\, \bm{Z}^{(0)}=\bm e_{d_{\ell'}})}
{{\P(\bm Z^{(n+1)}\neq \bm 0| \bm{Z}^{(0)}=\bm e_\ell)}}
\end{eqnarray*}
\end{lem}

\

\subsubsection*{Proof of Lemma \ref{GeigMultitype}}

%\begin{proof} 
The proof follows the same lines as the proof of Lemma 2.1 from \cite{G99}.  
Let $\mathcal T$ denote the tree of a branching process with the p.g.f of $Z$, and let $(A_i)_{1\leq i\leq |\bm z|}$ be measurable subsets of the space of multi-type rooted planar trees with roots of type $d_i$, such that for $1\leq i\leq j-1$
$$
A_i\subseteq \big\{\mathcal T\,:\,\bm Z^{(n)}(\mathcal T)=\bm0,\,\bm Z^{(0)}(\mathcal T)=\bm e_{d_i}\big\},\quad A_j\subseteq \big\{\mathcal T\,:\,\bm Z^{(n)}(\mathcal T)\neq\bm0,\,\bm Z^{(0)}(\mathcal T)=\bm e_{d_j}\big\},
$$
and with no condition on $A_i$ for $i>j$.  
Since $\{T^{(i)}\in A_i\}_{\,1\leq i\leq j-1}$ together with $T^{(j)}\in A_j$, imply that $R_{n+1}=j$, we have
$$
\Big\{\mathop{\cap}\limits_{i=1}^{|\bm z|}\big\{T^{(i)}\in A_i\big\},\, \bm Z^{(1)}=\bm z,\, \bm d(\bm Z^{(1)})=\bm d\Big\}\subset\big\{R_{n+1}=j\big\}
$$
so\vspace{-3mm}
\begin{eqnarray*} 
\P \big(\{T^{(i)}\in A_i\}_{\,1\leq i\leq |\bm z|}\!\!\!\!\!\!&,&\!\!\!\!\!\!\, \bm Z^{(1)}=\bm z,\,\bm d(\bm Z^{(1)})=\bm d,\,R_{n+1}=j\,\big|\bm Z^{(0)}=\bm e_\ell\big)\\
&=&\P\big(\bm Z^{(1)}=\bm z\big|\bm Z^{(0)}=\bm e_\ell\big)\P\big(\bm d\left(\bm z\right)=\bm d\big)\prod_{i=1}^{|\bm z|}\P\big(T^{(i)}\in A_i\,\big|\, \bm Z^{(0)}(T^{(i)})=\bm e_{d_i}\big).
\end{eqnarray*}
We next prove that for $1\leq i\leq j-1$
$$
\big(T^{(i)}\,\big|\,\bm Z^{(1)}=\bm z,\,\bm d(\bm Z^{(1)})=\bm d,\,\bm Z^{(0)}=\bm e_\ell, \, R_{n+1}=j
)\,\,\mathop{=}^d\,\,(\mathcal T\,\big|\,  \bm Z^{(n)}(\mathcal T)=\bm 0,\,\bm Z^{(0)}(\mathcal T)=\bm e_{d_i}\big).
$$
by showing that for every measurable subset $\hat{A}_i\subset \{\mathcal T\,|\, \bm Z^{(n)}(\mathcal T) =\bm 0, \, \bm Z^{(0)}(\mathcal T)=\bm e_{d_i}\big\}$, we have
$$
\P\big(T^{(i)}\in \hat{A}_i\big|\,\bm Z^{(1)}=\bm z,\,\bm d(\bm Z^{(1)})=\bm d,\,\bm Z^{(0)}=\bm e_\ell,  R_{n+1}=j \big)
=\P\big(\mathcal T\in \hat{A}_i\big| \bm Z^{(n)}(\mathcal T)=\bm 0,\,\bm Z^{(0)}(\mathcal T)=\bm e_{d_i}\big).
$$
The left hand side of the above equality can be rewritten as
\begin{equation}\label{prob1}
\frac{\P\big(T^{(i)}\in \hat{A}_i,\,\bm Z^{(1)}=\bm z,\,\bm d(\bm Z^{(1)})=\bm d,\,R_{n+1}=j\,\big|\,\bm Z^{(0)}=\bm e_\ell\big)}{\P\big(\bm Z^{(1)}=\bm z,\,\bm d(\bm Z^{(1)})=\bm d,\,R_{n+1}=j,\, \big|\,\bm Z^{(0)}=\bm e_\ell\big)}
\end{equation}
Using shorthand notation for events $E_i:= \big\{\mathcal T\,:\,\bm Z^{(n)}(\mathcal T)=\bm0,\,\bm Z^{(0)}(\mathcal T)=\bm e_{d_i}\big\}$, for $1\leq i\leq j-1$,\; 
$E_j:= \big\{\mathcal T\,:\,\bm Z^{(n)}(\mathcal T)\neq\bm 0,\,\bm Z^{(0)}(\mathcal T)=\bm e_{d_j}\big\},$ and \; 
$E_i:= \big\{\mathcal T\,:\,\bm Z^{(0)}(\mathcal T)=\bm e_{d_i}\big\}$,  for $j<i\leq|\bm z|$,
the numerator of (\ref{prob1}) becomes
\begin{eqnarray*} 
&\P&\!\!\!\!\!\big(T^{(i)}\in \hat{A}_i,\,\{T^{(r)}\in E_r\}_{r\neq i},\,\bm Z^{(1)}=\bm z,\,\bm d(\bm Z^{(1)})=\bm d,\,R_{n+1}=j\,\big|\,\bm Z^{(0)}=\bm e_\ell\big)\\ 
&=&\!\!\P\big(\bm Z^{(1)}=\bm z\big|\bm Z^{(0)}=\bm e_\ell\big)\,\P\big(\bm d(\bm z)=\bm d\big)\,\P\big(\mathcal T\in \hat A_i\big|\bm Z^{(0)}(\mathcal T)=\bm e_{d_i}\big)\prod_{r\neq i}\P\big(\mathcal T\in E_r\big|\bm Z^{(0)}(\mathcal T)=\bm e_{d_r}\big),
\end{eqnarray*}
while the denominator is equal to 
\begin{eqnarray*}
\P\big(\bm Z^{(1)}=\bm z,\,\bm d(\bm Z^{(1)})\!\!\!\!&=&\!\!\!\!\bm d,\, R_{n+1}=j\,\big|\bm Z^{(0)}=\bm e_\ell\big)\\ 
&=&\!\!\! \P\big(\bm Z^{(1)}=\bm z,\,\bm d(\bm Z^{(1)})=\bm d,\,\{T^{(r)}\in E_r\}_{r=1,\ldots,|\bm z|},\,R_{n+1}=j\,\big|\bm Z^{(0)}=\bm e_\ell\big)\\
&=&\!\!\! \P\big(\bm Z^{(1)}=\bm z\,\big|\,\bm Z^{(0)}=\bm e_\ell\big)\,\P\big(\bm d(\bm z)=\bm d\big)\prod_{r=1}^{|\bm z|}\P\big(\mathcal T\in E_r\big|\bm Z^{(0)}(\mathcal T)=\bm e_{d_r}\big).
\end{eqnarray*}
Together the last two equalities show that  (\ref{prob1}) is equal to 
$$
\P\big(T^{(i)}\in \hat{A}_i\big| R_{n+1}=j,\,\bm Z^{(1)}=\bm z,\,\bm d(\bm Z^{(1)})=\bm d,\,\bm Z^{(0)}=\bm e_\ell\big)=
\frac{\P\big(\mathcal T\in \hat A_i\big|\bm Z^{(0)}(\mathcal T)=\bm e_{d_i}\big)}{\P\big(\mathcal T\in E_j\big|\bm Z^{(0)}(\mathcal T)=\bm e_{d_j}\big)}
%\P\big(\mathcal T\in \hat{A}_i\big| \bm Z^{(n)}(\mathcal T)=\bm 0,\,\bm Z^{(0)}(\mathcal T)=\bm e_{d_i}\big).$$
$$
which by definition of $E_j$ is then equal to the right hand side of the equation above (\ref{prob1}).
Similar reasoning goes for for $i=j$ and for $i>j$. For the proof of the second part of the Lemma, it is sufficient to condition on $\bm Z_1=\bm z$ and use independence, as in proof for the single type case in \cite{G99}.\;\;$\Box$
%\end{proof}

\

\subsubsection*{\bf Proof of Theorem~\ref{Din-thm}:}\label{ProofDin-thm}
The proof of Markov property of the process $(\bm D_i)_{i\ge 1}$ is similar to the proof of the Markovian nature of $(D_i)_{i\ge 1}$ in Theorem 2.1 of \cite{LP13}. The main difference is that, instead of only recording the number of offspring  $D_i(n)$ of the ancestor of $(0,i)$ in generation $-n$ that have surviving progeny embedded to the right of $(0,i)$, we now record the types of these offspring including the one that is on the lineage of $(0,i)$ as well.
This, however, does not change the fact that 
\begin{eqnarray*}A_i>n&\Leftrightarrow& \forall n'\le n,\; a_i(n')\neq a_{i+1}(n')\\
&\Leftrightarrow& \forall n'\le n,\; (-n', a_i(n'))\mbox{ has no surviving progeny in }\{(0,i+1), (0,i+2),\dots\}\\
&\Leftrightarrow& \forall n'\le n,\; \|\bm D_i(n')\| =1
\end{eqnarray*}
so that $A_i$ is the level of the first term of the sequence $\bm D_i$ such that $D_i(n)=\|\bm D_i(n)\|-1\neq 0$.
In addition, at level $A_i$ we have the most recent common ancestor $a_i(A_i)=a_{i+1}(A_i)$ of individuals $(0,i)$ and $(0,i+1)$, whose offspring with surviving progeny embedded to the right of $(0,i+1)$ do not include the ancestor of $(0,i)$, which is recorded in $\bm D_i(A_i)_{[1]}$, but do include all the others. So,
$$\big(\bm D_{i+1}(A_i)_{[1]}, \bm D_{i+1}(A_i)_{[2]}, \bm D_{i+1}(A_i)_{[3]} \dots, \big)= \big(\bm D_{i}(A_i)_{[2]}, \bm D_{i}(A_i)_{[3]}, \dots\big).$$

At any level $n>A_i$ below the most recent common the ancestors of $(0,i)$ and $(0,i+1)$ are the same  
since $a_i(A_i)=a_{i+1}(A_i)$ implies $a_i(n)=a_{i+1}(n)$, so $$\forall n>A_i,\; \bm D_i(n)=\bm D_{i+1}(n).$$

For levels $n<A_i$ above the most recent common ancestor, note that the subtrees descending from different  surviving offspring of $(-A_i, a_i(A_i))$ are independent copies of multi-type branching processes whose initial individuals are of types $\bm D_i(A_i)_{[1]}, \bm D_i(A_i)_{[2]}, \dots$ and which are conditioned to survive for at least $n':=A_i-1$ generations. In particular, the subtree containing the lineage of $(0,i+1)$ above $(-A_i, a_i(A_i))$ is independent of the subtree whose lineage is recorded in $\{\bm D_i(n), n< A_i\}$ and is initiated by an individual of type $\ell:=\bm D_i(A_i)_{[2]}$. By definition $\big(\bm D_{i+1}(n), 1\le n<A_i\big)$ records the survivor types  (and their siblings) along the left most ancestral lineage of $(0,i+1)$ above the level $A_i$. The distribution of this sequence of type vectors for a multi-type branching process with initial individual of type $\ell$ conditioned to survive at least $n'$ generations is distributed as the sequence of type vectors $\bm \eta_{n',\ell}$. So,
$$\big(\bm D_{i+1}(1), \dots, \bm D_{i+1}(A_i-1)\big) \mathop{=}^d \bm{\eta}_{_{A_i-1,\bm{D}_{i}(A_i)_{[2]}}}\;\Leftrightarrow\;
\forall 1\le n<A_i,\; \bm D_{i+1}(n)\mathop{=}^d \bm{\eta}_{_{A_i-1,\bm{D}_{i}(A_i)_{[2]}}}(n).$$

As in the single type case, the sequence $\bm D_{i+1}=(\bm D_{i+1}(n), n\ge 1)$ depends only on $\bm D_i$ and not on $\bm D_{i'}$ for $i'<i$; and its transition law is determined by values of $(\bm D_i(n), \; n\ge A_i)$ and an independent random variable $\eta_{n', \ell}$ with $n'=A_i-1$ and $\ell=\bm D_i(A_i)_{[2]}$.\;\;$\Box$\\

\subsubsection*{\bf Proof of Proposition~\ref{DistribAi}:}\label{ProofDistribAi}
Observe that $\{A_1\neq 1,\dots, A_1\neq n\}$ iff all the ancestors $(-n',a_1(n'))$ of $(0,1)$ in generations $-n', \,1\le n'\le n$ have exactly one offspring with surviving progeny. 
By Lemma~\ref{GeigMultitype} when types of ancestral individuals are known the events of having exactly one offspring surviving progeny are independent across different generations. Let $\bm a_{[n']}=\bm A_{0[n']}$ denote the type of the ancestor $(-n', a_1(n'))$ of $(0,1)$ in generation $-n'$. These events can be expressed in terms of the random variable $\bm \eta_{n,\bm a_{[n]}}$ and in terms of the random variables  $\bm \zeta'_{n',\bm a_{[n']}},\, 1\le n'\le n$ as
\begin{eqnarray*}
\P\big(A_1>n, \bm A_{0|n-1}=\bm a_{|n-1}|\,\bm A_{0[n]}=\bm a_{[n]}\big)&=&\P\big(\bm \eta_{n,\bm a_{[n]}}=(\{\bm a_{[0]}\},\{\bm a_{[1]}\},\dots,\{\bm a_{[n-1]}\})\big)\\
&=&\P\big(\forall 1\le n'\le n\!:\, \bm \zeta'_{n',\bm a_{[n']}, \bm a_{[n'-1]}}=1,\; \bm \zeta'_{n',\bm a_{[n']}, \ell'}=0\; \forall \ell'\neq \bm a_{[n'-1]}\big)\\
&=&\prod_{n'=1}^n \P\big(\bm \zeta'_{n',\bm a_{[n']}, \bm a_{[n'-1]}}=1,\; \bm \zeta'_{n',\bm a_{[n']}, \ell'}=0\; \forall \ell'\neq \bm a_{[n'-1]}\big)
\end{eqnarray*}
where we can write the above as a product because the subtrees descending from different offspring are independent.
For each product term we have 
\begin{eqnarray*}
\P\big(\bm \zeta'_{n',\bm a_{[n']}, \bm a_{[n'-1]}}=1,\; \bm \zeta'_{n',\bm a_{[n']}, \ell'}=0\; \!\!\!\!\!\!&\forall&\!\!\!\!\!\! \ell'\neq \bm a_{[n'-1]}\big)\\
&=&\frac{\P(\bm \zeta_{n',\bm a_{[n']}, \bm a_{[n'-1]}}=1,\; \bm \zeta_{n',\bm a_{[n']}, \ell'}=0\; \forall \ell'\neq \bm a_{[n'-1]})}{\P(\sum_{\ell'=1}^k\bm \zeta_{n',\bm a_{[n']}, \ell'}\neq 0)}
\end{eqnarray*}
Conditioning on the value of variable $\bm \xi_{\bm a_{[n']}}$ which, when Bernoulli sampled by the vector $\bm p_{n'-1}:=\bm 1 -\bm f^{(n'-1)}(0,\ldots,0)$ of survival probabilities of different types by generation $n'-1$,  gives the distribution of $\bm \zeta_{n',\bm a_{[n']}}$, we get for the numerator
\begin{eqnarray*}
\P\big(\bm \zeta_{n',\bm a_{[n']}, \bm a_{[n'-1]}}=1,\; \bm \zeta_{n',\bm a_{[n']}, \ell'}=0\; \!\!\!\!\!\!&\forall&\!\!\!\!\!\! \ell'\neq \bm a_{[n'-1]}\big)\\
&=&\!\!\! \E\Big(\P\big(\bm \zeta_{n',\bm a_{[n']}, \bm a_{[n'-1]}}=1,\; \bm \zeta_{n',\bm a_{[n']}, \ell'}=0\; \forall \ell'\neq \bm a_{[n'-1]}\,|\bm \xi_{\bm a_{[n']}}\big)\Big)\\
&=&\!\!\!\E\Big(\bm \xi_{\bm a_{[n']},\bm a_{[n'-1]}}\frac{p_{n'-1,\bm a_{[n'-1]}}}{(1-p_{n'-1,\bm a_{[n'-1]}})}\prod_{\ell'=1}^k(1-p_{n'-1,\ell'})^{\bm \xi_{\bm a_{[n']},\ell'}}\Big),\\
&=&\!\!\!p_{n'-1,\bm a_{[n'-1]}}\frac{\partial f_{\bm a_{[n']}}(\bm s)}{\partial s_{\bm a_{[n'-1]}}}\Big|_{\bm s=\bm 1-\bm p_{n'-1}},
\end{eqnarray*}
and for the denominator 
\begin{eqnarray*}
\P\big(\sum_{\ell'=1}^k\bm \zeta_{n',\bm a_{[n']}, \ell'}\ge 1\big)
\!\!\!&=&\!\!\!1-\E\Big(\P\big(\sum_{\ell'=1}^k\bm \zeta_{n',\bm a_{[n']}, \ell'}= 0\,|\, \bm \xi_{\bm a_{[n']}}\big)\Big)
=1-\E\Big(\prod_{\ell'=1}^k(1-p_{n'-1,\ell'})^{\bm \xi_{\bm a_{[n']},\ell'}}\Big)\\
\!\!\!&=&\!\!\!1- f_{\bm a_{[n']}}(\bm 1-\bm p_{n'-1}) = 1- f^{(n')}_{\bm a_{[n']}}(0,\ldots,0)=p_{n',\bm a_{[n']}}.
\end{eqnarray*}
Since survival probability to generation $0$ is  $p_{0,\bm a_{[0]}}=1$, we have
\begin{eqnarray*}
\P(A_{1}>n,\bm A_{0|n-1}=\bm a_{|n-1}|\,\bm A_{0[n]}=\bm a_{[n]})
&=&\prod_{n'=1}^n \Big(
\frac{\partial f_{\bm a_{[n']}}(\bm s)}{\partial s_{\bm a_{[n'-1]}}} \Big|_{\bm s=\bm 1-\bm p_{n'-1}}
\frac{p_{n'-1,\bm a_{[n'-1]}}}{p_{n',\bm a_{[n']}}}\Big)\\
&=&\frac{1}{p_{n,\bm a_{[n]}}}\prod_{n'=1}^n \Big(\frac{\partial f_{\bm a_{[n']}}(\bm s)}{\partial s_{a_{[n'-1]}}} \Big|_{\bm s=\bm 1-\bm p_{n'-1}}\Big)
\end{eqnarray*}
and note that for $n'=1$ the evaluation of the derivative is at $\bm s=\bm 1-\bm p_{0}=\bm 0$. \;\; $\Box$\\

\

\subsubsection*{\bf Proof of Corollary~\ref{DistribAiuncond}:}\label{ProofDistribAiuncond}
One could obtain an expression for $\P(A_1>n|\, \bm A_{0[n]}=\bm a_{[n]})$ by summing over all the possible values for $\bm A_{0[n']}$ for $1\le n'<n$ that start with $\bm A_{0[n]}=\bm a_{[n]}$. This, however, in practice is only reasonable in special case of offspring distribution. On the other hand, the fact that $\{A_1\neq 1,\dots, A_1\neq n\}$ iff the subtree of the ancestor $(-n,a_1(n))$ of $(0,1)$ in generations $-n$ has exactly one offspring with surviving progeny, directly implies that
$$
\P(A_{1}>n |\,\bm A_{0[n]}=\bm a_{[n]})=\P\big( \sum_{\ell=1}^kZ^{(n)}_{\ell}=1\,|\, \bm Z^{(n)}\neq \bm 0, \bm Z^{(0)}=\bm e_{\bm a_{[n]}}\big ).\;\;\Box
$$

\

\subsubsection*{\bf Proof of Proposition~\ref{DistribBi}:}\label{ProofDistribBi}
Observe that $\{ B_{\ell,1}\neq 1, \dots, B_{\ell,1}\neq n\}$ iff all the ancestors $(-n',a_1(n'))$ of $(0,1)$ in generations $-n',1\le n'\le n$ have exactly one descendant in the standing population that has type $\ell$. As before, let $\bm a_{[n']}=\bm A_{0[n']}$ denote the type of the ancestor $(-n', a_1(n'))$ of $(0,1)$ in generation $-n'$, and note that $\bm a_{[0]}=\ell$.

We need to introduce new random variables which will count the number of offspring with descendants of type $\ell$ in the standing population. 
If $\bm f^{(n'-1)}=( f_1^{(n'-1)},\dots, f_k^{(n'-1)})$ is the probability generating function of the $n'-1$ generation in a multi-type branching process initiated by  individuals of type $\{1,\dots,k\}$, then the probability that a multi-type process after $n'-1$ generations has no individuals of type $\ell$ is given by the vector $\bm f^{(n'-1)}(\bm {\hat e}_\ell)$, where  $\bm {\hat e}_{\ell}:=\bm 1-\bm e_{\ell}$. Let $\bm p_{(n'-1,\ell)}:=\bm 1-\bm f^{(n'-1)}(\bm {\hat e}_\ell)$ denote the probability of having at least one descendant of type $\ell$ after $n'-1$ generations,  that is, for each $\jmath'\in\{1,\dots,k\}$, we have $p_{(n'-1,\ell),\jmath'}= 1-f^{(n'-1)}_{\jmath'}(\bm {\hat e}_\ell)$.
For all different offspring types $\jmath'\in\{1,\dots,k\}$, let $\{\varepsilon_{(n',\ell),\jmath'}^1,\varepsilon_{(n',\ell),\jmath'}^2,\dots \}$ be independent sequences of  Bernoulli variables with parameter $\P(\varepsilon_{(n',\ell),\jmath'}^m=1)=p_{(n'-1,\ell),\jmath'}$. For an independent offspring variable $\bm \xi_{\jmath}$ with $\bm \xi_{\jmath,\jmath'}$ offspring of type $\jmath'$ define 
$$\bm \zeta_{(n',\ell),\jmath}:=\Big(\sum_{m=1}^{\bm \xi_{\jmath, 1}} \varepsilon^m_{(n',\ell),1},\ldots, \sum_{m=1}^{\bm \xi_{\jmath, k}} \varepsilon^m_{(n',\ell),k}\Big)$$ 
whose $\jmath'$ coordinate is denoted by $\bm \zeta_{(n',\ell),\jmath,\jmath'}$. Then, $\bm \zeta_{(n',\ell),\jmath}$ records the number of offspring (of different types), in the first generation of a multi-type branching process initiated by an individual of type $\jmath$, which have at least one descendant of type $\ell$ after $n'-1$ generations. 

%Since, for each offspring, the event of having at least one descendant of type $\ell$ in current generation is a subset of the event of surviving, we can define a coupling of random variables  $\bm \zeta_{(n',\ell),\jmath}$ and $\bm \zeta_{n',\jmath}$ via a coupling of all the sequences of  Bernoulli variables $\{\varepsilon_{(n',\ell),\jmath'}^1,\varepsilon_{(n',\ell),\jmath'}^2,\dots \}$ used to define $\bm \zeta_{(n',\ell),\jmath}$ with respective sequences of Bernoulli variables $\{\epsilon_{n',\jmath'}^1,\epsilon_{n',\jmath'}^2,\dots \}$ used to define $\bm \zeta_{n',\jmath}$: define the first set of sequences from the second set of sequences using independent subsampling with probabilities $p_{(n'-1,\ell),\jmath'}/ p_{n'-1,\jmath'}$.
Then the law of $\bm \zeta'_{(n',\ell),\jmath}$ representing the number of offspring (of different types) in generation $-(n'-1)$ of a type $\jmath$ ancestor from generation $-n'$ whose descendants contain an individual of type $\ell$ in the standing population, given that there is at least one, is given by: 
%$$\bm \zeta'_{(n',\ell),\jmath}\mathop{:=}\limits^d\big(\bm \zeta_{(n',\ell),\jmath}\big|\sum_{\jmath'=1}^k\bm \zeta_{n',\jmath,\jmath'}\ge 1\big)$$
$$\bm \zeta'_{(n',\ell),\jmath}\mathop{:=}\limits^d\big(\bm \zeta_{(n',\ell),\jmath}\big|\sum_{\jmath'=1}^k\bm \zeta_{(n',\ell),\jmath,\jmath'}\ge 1\big)$$

The event $\{B_{\ell,1}>n\}=\{B_{\ell,1}\neq 1, \dots, B_{\ell,1}\neq n\}$ can now be expressed in terms of the newly defined random variables $\bm \zeta'_{(n',\ell),\bm a_{[n']}}, 1\le n'\le n$  as
\begin{eqnarray*}
\P\big(B_{\ell,1}>n\!\!\!\!\!\!&,&\!\!\!\!\! \bm A_{0|n-1}=\bm a_{|n-1}|\,\bm A_{0[n]}=\bm a_{[n]}, \bm A_{0[0]}=\ell \big)\\
&=&\P\big(\forall 1\le n'\le n: \bm \zeta'_{(n',\ell),\bm a_{[n']},\bm a_{[n'-1]}}=1,\; \bm \zeta'_{(n',\ell),\bm a_{[n']},\jmath'}=0 \,\forall \jmath'\neq \bm a_{[n'-1]}\big)\\
&=&\prod_{n'=1}^n \P\big(\bm \zeta'_{(n',\ell),\bm a_{[n']},\bm a_{[n'-1]}}=1,\; \bm \zeta'_{(n',\ell),\bm a_{[n']},\jmath'}=0 \,\forall \jmath'\neq \bm a_{[n'-1]}\big)
\end{eqnarray*}
where the product form follows since, by Lemma~\ref{GeigMultitype}, subtrees of different offspring are independent. For each product term we have
\begin{eqnarray*}
\P\big(\bm \zeta'_{(n',\ell),\bm a_{[n']},\bm a_{[n'-1]}}=1,\; \bm \zeta'_{(n',\ell),\bm a_{[n']},\jmath'}=0 \!\!\!\!\!&\forall&\!\!\!\!\! \jmath'\neq \bm a_{[n'-1]}\big)\\
&=& \frac{\P(\bm \zeta_{(n',\ell),\bm a_{[n']},\bm a_{[n'-1]}}=1, \bm \zeta_{(n',\ell),\bm a_{[n']}, \jmath'}=0\, \forall \jmath'\neq \bm a_{[n'-1]})}{\P(\sum_{\jmath'=1}^k\bm \zeta_{(n',\ell),\jmath,\jmath'}\ge 1)}
\end{eqnarray*}
A similar calculation to the one in the proof of Proposition~\ref{DistribAi}, conditioning on $\bm \xi_{\bm a_{[n']}}$, gives the numerator to be 
\begin{eqnarray*}
\P(\bm \zeta_{(n',\ell),\bm a_{[n']},\bm a_{[n'-1]}}=1, \bm \zeta_{(n',\ell),\bm a_{[n']},\jmath'}=0\!\!\!\!\!&\forall&\!\!\!\!\!\! \jmath'\neq \bm a_{[n'-1]})\\
&=&\!\!\! \E\Big(\P\big(\bm \zeta_{(n',\ell),\bm a_{[n']},\bm a_{[n'-1]}}=1, \bm \zeta_{(n',\ell),\bm a_{[n']},\jmath'}=0\,\forall\jmath'\neq \bm a_{[n'-1]}|\,\bm \xi_{\bm a_{[n']}}\big) \Big)\\
&=&\!\!\!\E\Big(\bm \xi_{\bm a_{[n']},\bm a_{[n'-1]}}\frac{p_{(n'-1,\ell),\bm a_{[n'-1]}}}{(1-p_{(n'-1,\ell),\bm a_{[n'-1]}})}\prod_{\jmath'=1}^k(1-p_{(n'-1,\ell),\jmath'})^{\bm \xi_{\bm a_{[n']},\jmath'}}\Big)\\
&=&p_{(n'-1,\ell),\bm a_{[n'-1]}}\frac{\partial f_{\bm a_{[n']}}(\bm s)}{\partial s_{\bm a_{[n'-1]}}}\Big|_{\bm s=\bm 1-\bm p_{(n'-1,\ell)}},
\end{eqnarray*}
while the denominator is calculated in the same way and equals
\begin{eqnarray*}
\P\big(\sum_{\ell'=1}^k\bm \zeta_{(n',\ell),\bm a_{[n']}, \ell'}\ge 1\big)
%&=&1-\P\big(\sum_{\ell'=1}^k\bm \zeta_{(n',\ell),\bm a_{[n']}, \ell'}= 0\big)\\
%\!\!\!&=&\!\!\!1-\E\Big(\P\big(\sum_{\ell'=1}^k\bm \zeta_{(n',\ell),\bm a_{[n']}, \ell'}= 0\,|\, \bm \xi_{\bm a_{[n']}}\big)\Big)
&=&1-\E\Big(\prod_{\ell'=1}^k(1-p_{(n'-1,\ell),\ell'})^{\bm \xi_{\bm a_{[n']},\ell'}}\Big)\\
%\!\!\!&=&\!\!\!
&=&1- f_{\bm a_{[n']}}(\bm 1-\bm p_{(n'-1,\ell)}) = 1- f^{(n')}_{\bm a_{[n']}}(\hat{\bm e}_{\ell})=p_{(n',\ell),\bm a_{[n']}}.
\end{eqnarray*}
%while the denominator is the same as in the product terms in the formula of Proposition~\ref{DistribAi}. 
%Note that now the factor $p_{(n'-1,\ell),\bm a_{[n'-1]}}$ in the numerator of the $n'$-th term and $p_{n'-1,\bm a_{[n'-1]}}$ in the denominator of the $(n'-1)$-th term no longer cancel, 
the cross terms of probabilities of a lineage with descendants of type $\ell$ cancel,  and 
$$\P(B_{\ell,1}>n, \bm A_{0|n-1}=\bm a_{|n-1}|\,\bm A_{0[n]}=\bm a_{[n]},\bm A_{0[0]}=\ell)= \frac{1}{p_{(n,\ell), \bm a_{[n]}}}
\prod_{n'=1}^n \frac{\partial f_{\bm a_{[n']}}(\bm s)}{\partial s_{\bm a_{[n'-1]}}}\bigg|_{\bm s=\bm 1-\bm p_{(n'-1,\ell)}}.\;\;\Box$$
%\end{proof}

\

\subsubsection*{\bf Proof of Corollary~\ref{DistribBiuncond}:}\label{ProofDistribBiuncond}
Again, an expression for $\P(B_{\ell,1}>n|\, \bm A_{0[n]}=\bm a_{[n]}, \bm A_{0[0]}=\ell)$ can be obtained by summing over all the possible values for $\bm A_{0[n']}$ for $1\le n'<n$ that start with $\bm A_{0[n]}=\bm a_{[n]}$, as can be seen in case of the LF  branching process. For a short expression we can use the fact that $\{B_{\ell,1}> n\}$ iff the subtree of the ancestor $(-n,a_1(n))$ of $(0,1)$ in generations $-n$ has exactly one descendant of type $\ell$ after $n$ generations, so
$$
\P(B_{\ell,1}>n |\,\bm A_{0[n]}=\bm a_{[n]},\bm A_{0[0]}=\ell)=\P\big( Z^{(n)}_{\ell}=1\,|\, Z^{(n)}_\ell\neq \bm 0, \bm Z^{(0)}=\bm e_{\bm a_{[n]}}\big ).\;\;\Box
$$

\

\subsection{Proof of results for LF multi-type branching process}
\subsubsection*{\bf Proof of Proposition~\ref{DistribAiLFCase2}:}\label{ProofDistribAiLFCase2}
% in conjunction with our
We start from our formula (\ref{A1law}) from Proposition~\ref{DistribAi},
\begin{eqnarray*}
\P(A_{1}>n|\,\bm A_{0[n]}=\bm a_{[n]})=\sum_{\bm a_{[1]},\dots,\bm a_{[n-1]}}
\prod_{n'=1}^n \Big(\frac{p_{n'-1,\bm a_{[n'-1]}}}{p_{n',\bm a_{[n']}}} \frac{\partial f_{\bm a_{[n']}}(\bm s)}{\partial s_{a_{[n'-1]}}} \Big|_{\bm s=\bm 1-\bm p_{n'-1}}\Big)
\end{eqnarray*}
in which we perform the summation in a 'top-down' order, from possible values for $\bm a_{[0]}$ down to $\bm a_{[n-1]}$.
Notice that only the first term in the product depends on $\bm a_{[0]}$ and write
\begin{eqnarray}\label{sumlayers}
\P(A_{1}>n|\,\bm A_{0[n]}=\bm a_{[n]})=\sum_{\bm a_{[1]},\dots,\bm a_{[n-1]}}\!
\Big(\sum_{\bm a_{[0]}=1}^k \frac{p_{0,\bm a_{[0]}}}{p_{1,\bm a_{[1]}}}\frac{\partial f_{\bm a_{[1]}}(\bm s)}{\partial s_{a_{[0]}}} \Big|_{\bm s=\bm 1-\bm p_{0}}\Big)
\prod_{n'=2}^n\frac{p_{n'-1,\bm a_{[n'-1]}}}{p_{n',\bm a_{[n']}}}\frac{\partial f_{\bm a_{[n']}}(\bm s)}{\partial s_{a_{[n'-1]}}} \Big|_{\bm s=\bm 1-\bm p_{n'-1}}
\end{eqnarray}
where the first summation is a function of $\bm a_{[1]}$ only, and equals
$$
c_1(\bm a_{[1]})=\sum_{\bm a_{[0]}=1}^k \frac{p_{0,\bm a_{[0]}}}{p_{1,\bm a_{[1]}}}\frac{\partial f_{\bm a_{[1]}}(\bm s)}{\partial s_{a_{[0]}}} \Big|_{\bm s=\bm 1-\bm p_{0}}= \frac{\P(\sum_{\ell'=1}^k\bm \zeta_{1,\bm a_{[1]},\ell'}=1)}{\P(\sum_{\ell'=1}^k\bm \zeta_{1,\bm a_{[1]}, \ell'}\ge 1)} 
$$
In general this value depends on $\bm a_{[1]}$, and the next `layer' of summation in (\ref{sumlayers}) is more complicated:
$$\sum_{\bm a_{[1]}=1}^k c_1(\bm a_{[1]})\frac{p_{1,\bm a_{[1]}}}{p_{2,\bm a_{[2]}}}\frac{\partial f_{\bm a_{[2]}}(\bm s)}{\partial s_{a_{[1]}}} \Big|_{\bm s=\bm 1-\bm p_{1}}$$
However, if $c_1(\bm a_{[1]})\equiv c_1$ for any allowed value of $\bm a_{[1]}$, then the next summation is simply
$$c_1\sum_{\bm a_{[1]}=1}^k \frac{p_{1,\bm a_{[1]}}}{p_{2,\bm a_{[2]}}}\frac{\partial f_{\bm a_{[2]}}(\bm s)}{\partial s_{a_{[1]}}} \Big|_{\bm s=\bm 1-\bm p_{1}}=c_1\cdot c_2(\bm a_{[2]})$$
If for all $1\le n'\le n$ we have $c_{n'}(\bm a_{[n']})\equiv c_{n'}$ for any allowed value of $\bm a_{[n']}$ then (\ref{sumlayers}) becomes simply\; $\P(A_1>n|\,\bm A_{0[n]}=\bm a_{[n]})=c_1\cdots c_n$,\; with
\begin{equation}\label{c}
c_{n'}(\bm a_{[n']})\equiv c_{n'}=\frac{1}{p_{n',\bm a_{[n']}}}\sum_{\bm a_{[n'-1]}=1}^k {p_{n'-1,\bm a_{[n'-1]}}}\frac{\partial f_{\bm a_{[n']}}(\bm s)}{\partial s_{a_{[n'-1]}}} \Big|_{\bm s=\bm 1-\bm p_{n'-1}}
\end{equation}
It also immediately follows that $\P(A_1>n|\, \bm A_{0[n]}=\bm a_{[n]})=\P(A_1>n)$.

\vspace{2mm}As we next show, this is precisely the case for multi-type LF branching processes.
The vector of extinction probabilities after $n'-1$ generations is $\bm p_{n'-1}=\! (1-h_{10}^{(n'-1)}\!,\dots,\!1-h_{k0}^{(n'-1)})$, while the derivative of the probability generating function evaluated at $\bm 1-\bm p_{n'-1}=\!(h_{10}^{(n'-1)}\!,\dots,\!h_{k0}^{(n'-1)})$ is 
\begin{eqnarray*}
\frac{\partial f_{\bm a_{[n']}}(\bm s)}{\partial s_{\bm a_{[n'-1]}}}\Big|_{\bm s=\bm 1-\bm p_{n'-1}}
\!\!\!\!=\frac{1}{1+m-m\sum_{\ell'=1}^k g_{\ell'} h_{\ell' 0}^{(n'-1)}}\Bigg(h_{\bm a_{[n']}0}+mg_{\bm a_{[n'-1]}}\frac{\sum_{\ell'=1}^k h_{\bm a_{[n']}\ell'} h_{\ell' 0}^{(n'-1)}}{1+m-m\sum_{\ell'=1}^k g_{\ell'} h_{\ell' 0}^{(n'-1)}}\Bigg)
\end{eqnarray*}
Let $U_{n'}=1+m-m\sum_{\ell'=1}^k g_{\ell'} h_{\ell' 0}^{(n'-1)}$ be the sum from the denominator, $V_{\bm a_{[n']}, n'}=\sum_{\ell'=1}^k h_{\bm a_{[n']}\ell'} h_{\ell' 0}^{(n'-1)}$ the sum from the numerator inside the bracket above. Then, we have 
\begin{eqnarray*}
c_{n'}(\bm a_{[n']})&=&\frac{1}{p_{n',\bm a_{[n']}}}\,\frac{\sum\limits_{\bm a_{[n'-1]}=1}^k (1-h^{(n'-1)}_{\bm a_{[n'-1]} 0})(h_{\bm a_{[n']}\bm a_{[n'-1]}}U_{n'}+mg_{\bm a_{[n'-1]}}V_{\bm a_{[n']},n'})}{U_{n'}^2}\\
&=&\frac{1}{p_{n',\bm a_{[n']}}}\,\frac{mV_{\bm a_{[n']},n'} (1-\!\!\!\sum\limits_{\bm a_{[n'-1]}=1}^k h^{(n'-1)}_{\bm a_{[n'-1]} 0}g_{\bm a_{[n'-1]}})+ U_{n'}(\!\!\!\sum\limits_{\bm a_{[n'-1]}=1}^k h_{\bm a_{[n']}\bm a_{[n'-1]}}-V_{\bm a_{[n']},n'})}{U_{n'}^2}
\end{eqnarray*}
Also note that 
$$
\frac{V_{\bm a_{[n']}, n'}}{U_{n'}}=f_{\bm a_{[n']}}(h_{10}^{(n'-1)}\!,\dots,\!h_{k0}^{(n'-1)}) -h_{\bm a_{[n']}0}=f_{\bm a_{[n']}}(\bm s)\big|_{\bm s=\bm 1-\bm p_{n'-1}} -h_{\bm a_{[n']}0}
$$
so, by the iterative property of the extinction probability, 
$$
p_{n',\bm a_{[n']}}=1-f_{\bm a_{[n']}}(\bm s)\big|_{\bm s=\bm 1-\bm p_{n'-1}}=1-h_{\bm a_{[n']}0}-\frac{V_{\bm a_{[n']}, n'}}{U_{n'}}=\frac{U_{n'}(1-h_{\bm a_{[n']}0})-V_{\bm a_{[n']},n'}}{U_{n'}}
$$
and after some arithmetic we get 
\begin{eqnarray*}
c_{n'}(\bm a_{[n']})&=&\frac{mV_{\bm a_{[n']},n'} (1-\!\!\!\sum\limits_{\bm a_{[n'-1]}=1}^k h^{(n'-1)}_{\bm a_{[n'-1]} 0}g_{\bm a_{[n'-1]}})+ U_{n'}(\!\!\!\sum\limits_{\bm a_{[n'-1]}=1}^k h_{\bm a_{[n']}\bm a_{[n'-1]}}-V_{\bm a_{[n']},n'})}{U_{n'}^2(1-h_{\bm a_{[n']}0})-U_{n'}V_{\bm a_{[n']},n'}}\\
&=&\frac{V_{\bm a_{[n']},n'} (U_{n'}-1)+ U_{n'}(1-h_{\bm a_{[n']}0}-V_{\bm a_{[n']},n'})}{U_{n'}^2(1-h_{\bm a_{[n']}0})-U_{n'}V_{\bm a_{[n']},n'}}=\frac{1}{U_{n'}}
\end{eqnarray*}
Hence, the formula (\ref{c}) for $c_{n'}(\bm a_{[n']})\equiv c_{n'}=(1+m-m\sum_{\ell'=1}^k g_{\ell'} h_{\ell' 0}^{(n'-1)})^{-1}$  is independent of $\bm a_{[n']}$ as claimed, and consequently
$$\P(A_1>n)=\P(A_1>n|\, \bm A_{0[n]}=\bm a_{[n]})=\prod_{n'=1}^n (1+m-m\sum_{\ell'=1}^k g_{\ell'} h_{\ell' 0}^{(n'-1)})^{-1}$$. 

\vspace{2mm}Independence of coalescence times $(A_i, i\ge 1)$ follows from the fact that in a multi-type LF branching process, all offspring other than the first one (which according to our current convention is the left-most one) are independent of the type of the parent, and have a Multivariate-Geometric distribution with mean $m^{(n)}$ and type distribution given by $\bm g^{(n)}$ (whose formula is given in  (\ref{LFparam})). 
This fact was also used in \cite{Sag13} (see Sec 4.1) to establish the formula (\ref{LFparam}) using the jumping contour representation of the branching process and its nice Markovian structure. 

We can rewrite this formula using the parameters defined in (\ref{LFparam}) from Theorem~\ref{thm-LF} according to which $$\bm H^{(n)}\bm =\bm M^n\bm -\frac{m^{(n)}}{1+m^{(n)}}\bm M^n \bm 1^{\t} \bm g^{(n)}$$
which when multiplied by $\bm g$ on the left and by $\bm 1^{\t}$ on the right becomes
$\bm g \bm H^{(n)}\bm 1^t=\frac{\bm g\bm M^n\bm 1^t}{1+m^{(n)}}$ or 
$\sum_{\ell=1}^kg_\ell h_{\ell 0}^{(n)}=1-\frac{\bm g \bm M^n\bm 1^t}{1+m^{(n)}}.$
Using this equality in the formula for $\P(A_1>n)$ we get
\begin{eqnarray*}
\P(A_1>n)\!\!\!\!&=&\!\!\!\!
\prod_{n'=1}^n \frac{1}{1+m-m\sum_{i=1}^k g_ih_{i0}^{(n'-1)}}=
\prod_{n'=1}^n \frac{1}{1+m\frac{\bm g \bm M^{n'-1}\bm 1^t}{1+m^{(n'-1)}}}\\
\!\!\!\!&=&\!\!\!\!\prod_{n'=1}^n \frac{1+m^{(n'-1)}}{1+m^{(n'-1)}+m\bm g \bm M^{n'-1}\bm 1^{\t}}
=\prod_{n'=1}^n \frac{1+m^{(n'-1)}}{1+m^{(n')}}
=\frac{1}{1+m^{(n)}}.
\end{eqnarray*}
\vspace{1mm}because $m^{(n'-1)}+m\bm g \bm M^{n'-1}\bm 1^{\t}=m\bm g(\bm I+\bm M+\cdots+\bm M^{n'-2})\bm 1^{\t}+m\bm g \bm M^{n'-1}\bm 1^{\t}=m^{(n')}$, and $m^{(0)}=0$.
Another way to see why this holds is to notice that %$\{A_1>n\}$ iff there is exactly one surviving individual in generation $n$ for a multi-type branching process started by the ancestral individual $(-n,a_1(n))$. Hence, 
$$\P(A_1>n|\,  \bm A_{0[n]}=\bm a_{[n]})=\P(\sum_{\ell=1}^kZ^{(n)}_{\ell}=1\,|\, \bm Z^{(n)}\neq \bm 0, \bm Z^{(0)}= \bm e_{\bm a_{[n]}})$$ which, as a result of Theorem~\ref{thm-LF} is simply the probability the Geometric variable with mean $m^{(n)}$ is 0, which is equal to $1/(1+m^{(n)})$ regardless of the type of the initial individual $\bm Z^{(0)}$. Although the latter approach is much shorter, we thought it would be instructive to show the agreement with the formula for the joint law of $A_1$ and $\bm A_{0|n-1}$ via the summation approach. \;$\Box$\\\\

\subsubsection*{\bf Proof of Proposition~\ref{DistribBiLFCase2}:}\label{ProofDistribBiLFCase2}
As in the proof of Proposition~\ref{DistribAiLFCase2} we show that in our formula (\ref{B1law}) from Proposition~\ref{DistribBi},
$$\P(B_{\ell,1}>n|\,\bm A_{0[n]}=\bm a_{[n]},\bm A_{0[0]}=\ell)=\sum_{\bm a_{[1]},\dots,\bm a_{[n-1]}}\!\prod_{n'=1}^n \Big(\frac{p_{(n'-1,\ell),\bm a_{[n'-1]}}}{p_{(n',\ell), \bm a_{[n']}}}\frac{\partial f_{\bm a_{[n']}}(\bm s)}{\partial s_{\bm a_{[n'-1]}}}\Big|_{\bm s=\bm 1-\bm p_{(n'-1,\ell)}}\Big)$$
for each $1\le n'\le n-1$ the sums 
\begin{equation}\label{tildec}
\tilde c_{n'}(\bm a_{[n']})\equiv \tilde c_{n'}= \frac{1}{p_{(n',\ell), \bm a_{[n']}}}\sum_{\bm a_{[n'-1]}=1}^k {p_{(n'-1,\ell),\bm a_{[n'-1]}}}\frac{\partial f_{\bm a_{[n']}}(\bm s)}{\partial s_{\bm a_{[n'-1]}}}\Big|_{\bm s=\bm 1-\bm p_{(n'-1,\ell)}}
\end{equation}
are independent of the value of $\bm a_{[n']}$. 
The only difference in the formula (\ref{tildec}) is that the fractional factors and the evaluation of the derivative now use $\bm 1-\bm p_{(n'-1,\ell)}:=\bm f^{(n'-1)}(\bm {\hat e}_{\ell'})$,  with $\bm {\hat e}_{\ell'}=\bm 1-\bm e_{\ell'}$, instead of using $\bm 1-\bm p_{n'}:=\bm f^{(n'-1)}(\bm 0)$ as they did in the formula (\ref{c}) for $c_{n'}$. In order to reuse the calculations we did for  $c_{n'}$, we define the following analogous notation:    
$$
\tilde h^{(n'-1)}_{\bm a_{[n'-1]}0}:=1-p_{(n'-1,\ell),\bm a_{[n'-1]}}=h^{(n'-1)}_{\bm a_{[n'-1]}0}+\frac{1-h^{(n'-1)}_{\bm a_{[n'-1]}0}-h^{(n'-1)}_{\bm a_{[n'-1]}\ell}
%\sum_{\jmath\neq \ell} h^{(n'-1)}_{\bm a_{[n'-1]}\jmath}
}{1+m^{(n'-1)}g^{(n'-1)}_{\ell}},
$$
$$\tilde U_{n'}:=1+m-m\sum_{\ell'=1}^k g_{\ell'} \tilde h_{\ell' 0}^{(n'-1)},\quad \tilde V_{\bm a_{[n']}, n'}=\sum_{\ell'=1}^k h_{\bm a_{[n']}\ell'} \tilde h_{\ell' 0}^{(n'-1)}.$$ Using $\tilde h^{(n'-1)}_{\bm a_{[n'-1]}0}$, $\tilde U_{n'}$ and $\tilde V_{\bm a_{[n']},n'}$ in place of $h^{(n'-1)}_{\bm a_{[n'-1]}0}$, $U_{n'}$ and $V_{\bm a_{[n']},n'}$ respectively, the exact same arithmetic follows through as in the proof of Proposition~\ref{DistribAiLFCase2} and (\ref{tildec}) becomes
$$
\tilde c_{n'}(\bm a_{[n']})= \frac{1}{\tilde U_{n'}}= \frac{1}{1+m-m\sum_{\ell'=1}^k g_{\ell'} \tilde h_{\ell' 0}^{(n'-1)}}
$$
and consequently
$$\P(B_{\ell,1}>n|\, \bm A_{0[0]}=\ell)=\P(B_{\ell,1}>n|\, \bm A_{0[n]}=\bm a_{[n]}, \bm A_{0[0]}=\ell)=\prod_{n'=1}^n (1+m-m\sum_{\ell'=1}^k g_{\ell'} \tilde h_{\ell' 0}^{(n'-1)})^{-1}$$
where for each $\ell'\in\{1,\dots,k\}$
$$ \tilde h_{\ell' 0}^{(n'-1)}=h^{(n'-1)}_{\ell'0}+\frac{1-h^{(n'-1)}_{\ell'0}-h^{(n'-1)}_{\ell'\ell}}{1+m^{(n'-1)}g^{(n'-1)}_{\ell}}.
$$
Independence of coalescence times $(B_{\ell,i})_{i\ge 1}$ follows by the same arguments as for $(A_i)_{ i\ge 1}$.

We can rewrite this formula in a similar way as before by noting that the first product term in (\ref{B1lawLFCase2}) is equal to ${1}/{(1+mg_{\ell})}$, and using many arithmetic steps established by the relationship of parameters in (\ref{LFparam}) the rest of the terms for $n'>1$  can be shown to be equal to 
$$ \frac{(1+m^{(n'-1)})(1+m^{(n'-1)}g_\ell^{(n'-1)})}{(1+m^{(n'-1)})(1+m\bm g\bm H^{(n'-1)}\bm e_\ell)+m^{(n'-1)}g_\ell^{(n'-1)}(1+m^{(n')})}=\frac{1+m^{(n'-1)}g_\ell^{(n'-1)}}{1+m^{(n')}g_\ell^{(n')}}$$
hence
$$\P(B_{\ell,1}>n|\, \bm A_{0[0]}=\ell)=\frac{1}{1+mg_{\ell}}\prod_{n'=2}^n \frac{1+m^{(n'-1)}g_\ell^{(n'-1)}}{1+m^{(n')}g_\ell^{(n')}}=\frac{1}{1+m^{(n)}g_\ell^{(n)}}.$$\\
Another way to see this expression is to notice that
 $$\P(B_{\ell,1}>n|\,  \bm A_{0[n]}=\bm a_{[n]},  \bm A_{0[0]}=\ell)=\P(Z^{(n)}_{\ell}=1\,|\, Z^{(n)}_{\ell}\neq \bm 0, \bm Z^{(0)}= \bm e_{\bm a_{[n]}})$$ which is equal to ${\partial_{s_\ell} f^{(n)}_{\bm a_{[n]}}(\hat{\bm e}_{\ell})}/{(1-f^{(n)}(\hat{\bm e}_{\ell}))}$, and as a result of Theorem~\ref{thm-LF} and some simple arithmetic can be shown to be  equal to $1/(1+m^{(n)}g^{(n)}_\ell)$ regardless of the type of the initial individual $\bm Z^{(0)}$.\;$\Box$\\

\subsubsection*{\bf Proof of Corollary~\ref{DistribLFCase1}}\label{ProofDistribLFCase1}
One approach is to use equivalence of offspring laws for different parents. Since types do not affect the reproduction law, the ancestral tree shape can be decoupled from the individual types. In other words,  we can first construct the branching process using the single-type LG offspring distribution with parameters $(h_0,m)$, and subsequently assign types to all individuals independently according to probabilities $\bm g$. As the only factor affecting the coalescent times of the standing population is the offspring number of each individual, their law is the same as in the case of the associated single-type LF coalescent times.

According to Proposition 5.1 from \cite{LP13} for the single-type LF case, the coalescent times have distribution (in their notation $a\mapsto h_0, b\mapsto m/(1+m)$ and offspring mean $m\mapsto (1+m)(1-h_0)$)
$$\P(A_1>n)=\frac{m-h_0(1+m)}{m(1+m)^n(1-h_0)^n-h_0(1+m)},$$  if $(1+m)(1-h_0)\neq 1$, while if $(1+m)(1-h_0)=1$
%$$\P(A_1>n)=\frac{1-(1+m)(1-h_0)}{1-(1+m)(1-h_0)+m(1-(1+m)^n(1-h_0)^n)},\quad\mbox{ if } (1+m)(1-h_0)\neq 1$$ and 
$$ \P(A_1>n)=\frac{1-h_0}{1-h_0+nh_0}.$$

To see that this agrees with result (\ref{A1lawLFCase2}) note $\bm H=(1-h_0)\bm 1^{\t}\bm g$, $\bm M=\bm H +m\bm H\bm 1^{\t}\bm g$ implies
$$\bm M=(1-h_0)\bm g +m(1-h_0)\bm 1^{\t}\bm g\bm 1^{\t} \bm g=(1-h_0)(1+m)\bm 1^{\t}\bm g,\;\;\; \bm M^{n'}=(1-h_0)^{n'}(1+m)^{n'}\bm 1^{\t}\bm g$$
\begin{eqnarray*}m^{(n)}\!\!\!\!&=&\!\!\!\!%m\bm g(\bm I+\bm M+\cdots +\bm M^{n-1})\bm 1^{\t}=
m\bm g(\bm I+(1-h_0)(1+m)\bm 1^{\t}\bm g+\dots+(1-h_0)^{n-1}(1+m)^{n-1}\bm 1^{\t}\bm g)\bm 1^{\t}\\
\!\!\!\!&=&\!\!\!\! m\sum_{n'=0}^{n-1}(1-h_0)^{n'}(1+m)^{n'}=m\frac{1-(1-h_0)^{n}(1+m)^{n}}{1-(1-h_0)(1+m)},\quad \mbox{ if } (1-h_0)(1+m)\neq 1
\end{eqnarray*}
while if $(1-h_0)(1+m)=1$ then $\bm M^{n'}=\bm M, m^{(n)}=mn=nh_0/(1-h_0)$. Using this in (\ref{A1lawLFCase2}) the formula  $\P(A_1>n)=1/(1+m^{(n)})$ gives the same result as above.

\vspace{5mm} From the result for $A_1$ we can obtain the law of $B_{\ell,1}$ using its original definition as the maximum of all coalescence times until the first next individual in the current population whose type is $\ell$: $B_{\ell, 1}:=\max\{A_{\i_{\ell,1}},\ldots,A_{\i_{\ell,2}-1}\}$. The only reason why this calculation is simple is due to the decoupling of the branching tree and the individual types. Since, given the branching tree, all individuals are assigned types independently according to probabilities $\bm g$, the index $\i_{\ell,2}:=\min\{i'> \i_{\ell,1}\equiv 1: \bm A_{i'[0]}=\ell\}$ is such that $\i_{\ell,2}-1$ is a random variable with a shifted Geometric distribution with parameter $g_\ell$.  Conditioning on the value of $\i_{\ell,2}-1$, and using the fact that $(A_i)_{i\ge 1}$ is an i.i.d sequence, we get
\begin{eqnarray*}
\P(B_{\ell,1}\le n)&=&\E\Big(\P\big(\max\{A_1,\dots,A_{\i_{\ell,2}-1}\}\le n \big|\,\i_{\ell,2}-1\big)\Big)\\
&=&\sum_{i'-1=1}^{\infty} \P(A_1\le n)^{i'-1} (1-g_\ell)^{i'-2}g_\ell =\frac{g_\ell \P(A_1\le n)}{1-\P(A_1\le n) (1-g_\ell)}
\end{eqnarray*}
and
$$\P(B_{\ell,1}>n)=1-\P(B_{\ell,1}\le n)=\frac{1-\P(A_1\le n)}{1-\P(A_1\le n)(1-g_\ell)}.\;\Box$$

\

\subsection{Calculations for two-type LF branching process}
\subsubsection*{Proof of Proposition~\ref{2typecomp}}\label{Proof2typecomp}
Formulae (\ref{LFparam}) imply (after much arithmetic using Maple) that
$$m^{(n)}=m_{\s}^{(n)}=m_{\a}^{(n)}=m\sum_{n'=0}^{n-1}h_1^{n'}(m+1)^{n'}=\frac{m(h_1^n(1+m)^{n}-1)}{h_1(1+m)-1},$$ 
$$\bm g_{\s}^{(n)}=\Big(\big(g-\frac{1}{2}\big)G(2p-1)+\frac{1}{2}, -\big(g-\frac{1}{2}\big)G(2p-1)+\frac{1}{2}\Big), \quad\bm g_{\a}^{(n)}=(gG(p), -gG(p)+1),\vspace{2mm}$$
where $G$ is a rather complicated polynomial
$$G(x)=\frac{h_1^{n-1}(h_1(1+m)-1)}{h_1^n(1+m)^n-1}\bigg[\sum_{n'=0}^{n-2}\Big(\frac{h_1m+1}{h_1^{n-n'-1}}+\sum_{i=1}^{n-n'-2}(1+m)^ih_1^{i-n+n'+2}\Big)x^{n'}+x^{n-1}\bigg],$$
which satisfies $G(0)=0, G(1)=1$. %\textcolor{blue}{\it is $G(x)$ convex, could that be useful??}\\
From (\ref{A1lawLFCase2}) we have the same formulae for the distribution of coalescent times in the two cases: 
$$\P_{\s}(A_1>n)=\P_{\a}(A_1>n)=\big(1+m^{(n)}\big)^{-1},$$%\big(1+m\sum_{n'=0}^{n-1}h_1^{n'}(m+1)^{n'}\big)^{-1},$$
and from (\ref{B1lawLFCase2}) we get the following formulae for same-type coalescent times in the two cases:
$$\P_{\s}(B_{\ell,i}>n\,|\,\bm A_{0[0]}=\ell)=\big(1+m^{(n)}g_{\s\ell}^{(n)}\big)^{-1},%\big(1+m\sum_{{n'}=0}^{n-1}h_1^{n'}(m+1)^{n'}g_{\s\ell}^{(n)}\big)^{-1},
\; \mbox{ for }\; \ell\in\{\bm 1, \bm 2\}$$
%$$\P_{\s}(B_{\bm 2,i}>n\,|\, \cdot)=\big(1+m\sum_{{n'}=0}^{n-1}h_1^{n'}(m+1)^{n'}g_{\s2}^{(n)}\big)^{-1}$$
$$\P_{\a}(B_{\ell, i}>n\,|\, \bm A_{0[0]}=\ell)=\big(1+m^{(n)}g_{\a\ell}^{(n)}\big)^{-1},%\big(1+m\sum_{{n'}=0}^{n-1}h_1^{n'}(m+1)^{n'}g_{\a\ell}^{(n)}\big)^{-1},
\;\mbox{ for }\; \ell\in\{\bm 1, \bm 2\}$$
%$$\P_{\a}(B_{\bm 2, i}>n\,|\, \cdot)=\big(1+m\sum_{{n'}=0}^{n-1}h_1^{n'}(m+1)^jg_{\a2}^{(n)}\big)^{-1}$$
where the two coordinates of vectors $ \bm g_{\s}^{(n)}=(g_{\s1}^{(n)},g_{\s2}^{(n)})$ and $ \bm g_{\a}^{(n)}=(g_{\a1}^{(n)},g_{\a2}^{(n)})$ are given above.\\
We next prove that for $g\leq 1/2$ we have
$$1+m^{(n)}g_{\a1}^{(n)}\leq 1+m^{(n)}g_{\s1}^{(n)}\; \Leftrightarrow\; g_{\a1}^{(n)}\leq g_{\s1}^{(n)},\quad 1+m^{(n)}g_{\s2}^{(n)}\leq1+m^{(n)}g_{\a2}^{(n)} \; \Leftrightarrow\;g_{\s2}^{(n)}\leq g_{\a2}^{(n)}.$$
Both of these inequalities are equivalent to
$$gG(p)+\big(\frac{1}{2}-g\big)G(2p-1)-\frac{1}{2}\leq 0,$$
which holds since all multiplying coefficients of the polynomial $G(x)$ are nonnegative, so $G$ is increasing and both $G(p),G(2p-1)\le G(1)=1$. \\

\noindent For the last comparison we need to show that for $g\leq 1/2\leq p$
$$1+m^{(n)}g_{\s1}^{(n)}\leq1+m^{(n)}g_{\s2}^{(n)}\; \Leftrightarrow\; g_{\s1}^{(n)}\leq g_{\s2}^{(n)}$$ 
which is equivalent to
$$2\Big(\frac12-g\Big)G(2p-1)\geq0,$$
and holds as long as $p\geq 1/2$ so that the polynomial $G(x)$ is evaluated on $x\ge 0$.\;$\Box$
%\newpage

\

\

\end{document}